\documentclass[11pt,reqno,final]{amsart}
\usepackage{amsmath,amsthm,amssymb,verbatim,enumerate,ifthen,times}
\usepackage[mathscr]{eucal}
\usepackage{showkeys}
\usepackage[utf8]{inputenc}
\newcommand{\N}{\mathbb{N}}
\newcommand{\R}{\mathbb{R}}
\newcommand{\C}{\mathscr{C}}
\newcommand{\CM}{\mathscr{CM}}
\newcommand{\CP}{\mathscr{CP}}
\newcommand{\A}{\mathscr{A}}
\newcommand{\B}{\mathscr{B}}

\newcommand{\D}{\mathscr{D}}

\def\d{\,{\rm d}}

\newtheorem{theorem}{Theorem}
\newtheorem*{theorem*}{Theorem}
\def\Thm#1#2{\ifthenelse{\equal{#1}{*}}{\begin{theorem*}#2\end{theorem*}}
  {\begin{theorem}\label{T#1}#2\end{theorem}}}
\newtheorem{Atheorem}{Theorem}

\def\thm#1{Theorem~\ref{T#1}}

\newtheorem{proposition}[theorem]{Proposition}
\newtheorem*{proposition*}{Proposition}
\def\Prp#1#2{\ifthenelse{\equal{#1}{*}}{\begin{proposition*}#2\end{proposition*}}
             {\begin{proposition}\label{P#1}#2\end{proposition}}}

\newtheorem{corollary}[theorem]{Corollary}
\newtheorem*{corollary*}{Corollary}
\def\Cor#1#2{\ifthenelse{\equal{#1}{*}}{\begin{corollary*}#2\end{corollary*}}
             {\begin{corollary}\label{C#1}#2\end{corollary}}}
\def\cor#1{Corollary~\ref{C#1}}

\newtheorem{lemma}[theorem]{Lemma}
\newtheorem*{lemma*}{Lemma}
\def\Lem#1#2{\ifthenelse{\equal{#1}{*}}{\begin{lemma*}#2\end{lemma*}}
             {\begin{lemma}\label{L#1}#2\end{lemma}}}
\def\lem#1{Lemma~\ref{L#1}}

\newtheorem{example}[theorem]{Example}
\newtheorem*{example*}{Example}
\def\Exa#1#2{\ifthenelse{\equal{#1}{*}}{\begin{example*}\rm #2\end{example*}}
             {\begin{example}\label{Ex#1}\rm #2\end{example}}}

\newtheorem{problem}[theorem]{Problem}

\theoremstyle{definition}
\newtheorem{definition}[theorem]{Definition}

\newtheorem{remark}[theorem]{Remark}
\newtheorem*{remark*}{Remark}
\def\Rem#1#2{\ifthenelse{\equal{#1}{*}}{\begin{remark*}\rm #2\end{remark*}}
             {\begin{remark}\label{R#1}\rm #2\end{remark}}}

\newcommand{\eq}[1]{\eqref{E#1}}
\newcommand{\Eq}[2]{\ifthenelse{\equal{#1}{*}}
  {\begin{equation*}\begin{aligned}[]#2\end{aligned}\end{equation*}}
  {\begin{equation}\begin{aligned}[]\label{E#1}#2\end{aligned}\end{equation}}}
\newcommand{\DET}[1]{\begin{vmatrix}#1\end{vmatrix}}
\def\lege{\begin{array}{c}<\\[-2mm]=\\[-2mm]>\end{array}}
\long\def\comment#1{}

\begin{document}
\large

\date{\today}

\title[Comparison of generalized Bajraktarevi\'c means]
{On the local and global comparison of generalized Bajraktarevi\'c means}

\author[Zs. P\'ales]{Zsolt P\'ales}
\address{Institute of Mathematics, University of Debrecen,
H-4002 Debrecen, Pf.\ 400, Hungary}
\email{pales@science.unideb.hu}

\author[A. Zakaria]{Amr Zakaria}
\address{Department of Mathematics, Faculty of Education, Ain Shams University, Cairo 11341, Egypt}
\email{amr.zakaria@edu.asu.edu.eg}


\thanks{This research has been supported by the Hungarian
Scientific Research Fund (OKTA) Grants K-111651.}

\begin{abstract}
Given two continuous functions $f,g:I\to\R$ such that $g$ is positive and $f/g$ is strictly monotone, a
measurable space $(T,\A)$, a measurable family of $d$-variable means $m: I^d\times T\to I$, and a probability
measure $\mu$ on the measurable sets $\A$, the $d$-variable mean $M_{f,g,m;\mu}:I^d\to I$ is defined by
\Eq{*}{
   M_{f,g,m;\mu}(\pmb{x})
      :=\left(\frac{f}{g}\right)^{-1}\left(
            \frac{\int_T f\big(m(x_1,\dots,x_d,t)\big)\d\mu(t)}
                 {\int_T g\big(m(x_1,\dots,x_d,t)\big)\d\mu(t)}\right)
   \qquad(\pmb{x}=(x_1,\dots,x_d)\in I^d).
}
The aim of this paper is to study the local and global comparison problem of these means, i.e., to find 
conditions for the generating functions $(f,g)$ and $(h,k)$, for the families of means $m$ and $n$, and for 
the measures $\mu,\nu$ such that the comparison inequality
\Eq{*}{
   M_{f,g,m;\mu}(\pmb{x})\leq M_{h,k,n;\nu}(\pmb{x}) \qquad(\pmb{x}\in I^d)
}
be satisfied.
\end{abstract}
\maketitle

\section{Introduction}

In a recent paper \cite{LosPal08}, Losonczi and P\'ales investigated a general class of two-variable means 
given by the formula
\Eq{*}{
   M_{f,g;\mu}(x_1,x_2)
      :=\left(\frac{f}{g}\right)^{-1}\left(
            \frac{\int_0^1 f(tx_1+(1-t)x_2)\d\mu(t)}{\int_0^1 g(tx_1+(1-t)x_2)\d\mu(t)}\right)
   \qquad((x_1,x_2)\in I^2),
}
where $f,g:I\to\R$ are continuous functions such that $g$ is positive and $f/g$ is strictly monotone and 
$\mu$ is a probability measure on the Borel measurable subsets of $[0,1]$. This definition includes many 
former classical and important settings. In \cite{LosPal08} local and global comparison theorems 
(that provided necessary and in some cases also sufficient conditions) have been established for the 
comparison of two two-variable means from this general class. 

The purpose of this note is to extend the results of \cite{LosPal08} in several ways. In our approach we will 
use Chebyshev systems, measurable families of means and measures for the definition of a general class 
of $d$-variable means.  

Throughout this paper, the symbols $\N$, $\R$, and $\R_+$ will stand for the sets of natural,
real, and positive real numbers, respectively, and $I$ will always denote a nonempty open real interval. The
classes of continuous strictly monotone and continuous positive real-valued functions defined on $I$ will be
denoted by $\CM(I)$ and $\CP(I)$, respectively.

In the sequel, a function $M:I^d\to I$ is called a
\emph{$d$-variable mean} on $I$ if the following so-called mean value property
\Eq{1}{
  \min(x_1,\dots,x_d)\leq M(\pmb{x})\leq \max(x_1,\dots,x_d)   \qquad(\pmb{x}=(x_1,\dots,x_d)\in I^d)
}
holds. Also, if both of the inequalities in \eq{1} are strict for all $x_1,\dots,x_d\in I$ with $x_i\neq x_j$
for some $i\neq j$, then we say that $M$ is a \emph{strict mean} on $I$. The \emph{arithmetic} and
\emph{geometric} means are well known instances for strict means on $\R_{+}$. More generally, if $p$ is a
real number, then the $d$-variable \textit{H\"older mean} $H_p:\R_+^d\to\R$ is defined as
\Eq{*}{
  H_{p}(\pmb{x})
  :=\begin{cases}
    \left(\dfrac{x_1^p+\cdots+x^p_d}{d}\right)^{\frac{1}{p}} &\mbox{if } p\neq0\\[3mm]
    \sqrt[d]{x_1\cdots x_d} &\mbox{if } p=0
   \end{cases}
  \qquad\big(\pmb{x}=(x_1,\dots,x_d)\in\R_+^d\big).
}
Obviously, $H_1$ and $H_0$ equal the arithmetic and geometric mean, respectively. It is easy to see that
H\"older means are strict means. The $d$-variable minimum and maximum functions are instances for non-strict
means.

A classical generalization of H\"older means is the notion of \emph{$d$-variable quasi-arithmetic mean} (cf.\
\cite{HarLitPol34}), which is introduced as follows: For $f\in\CM(I)$ define
\Eq{A}{
  A_f(\pmb{x}):=f^{-1}\left(\frac{f(x_1)+\cdots+f(x_d)}{d}\right)
  \qquad\big(\pmb{x}=(x_1,\dots,x_d)\in I^d\big).
}
More generally, if $S_d$ denotes the $(d-1)$-dimensional simplex given by
\Eq{S}{
  S_d:=\{(t_1,\dots,t_d)\mid t_1,\dots,t_d\geq0,\,t_1+\dots+t_d=1\},
}
then we can also define
\Eq{At}{
  A_f(\pmb{x},\pmb{t}):=f^{-1}\big(t_1f(x_1)+\cdots+t_df(x_d)\big)
  \qquad\big(\pmb{x}=(x_1,\dots,x_d)\in I^d,\, \pmb{t}=(t_1,\dots,t_d)\in S_d\big),
}
which is called the \emph{weighted $d$-variable quasi-arithmetic mean on $I$.}

In this paper, we consider a much more general class of means. For their definition, we recall the notion of
Chebyshev system. Let $f,g:I\to\R$ be continuous function. We say that the pair $(f,g)$ forms a
\textit{(two-dimensional) Chebyshev system on $I$} if, for any distinct elements $x,y$ of $I$, the
determinant
\Eq{*}{
  \D_{f,g}(x,y):=\DET{f(x) & f(y) \\ g(x) & g(y)}\qquad(x,y\in I)
}
is different from zero. If, for $x<y$, this determinant is positive, then $(f,g)$ is called a \emph{positive
system}, otherwise we call $(f,g)$ a \emph{negative system}. Due to the connectedness of the triangle
$\{(x,y)\mid x<y,\,x,y\in I\}$, it follows that every Chebyshev system is either positive or negative.
Obviously, if $(f,g)$ is a positive Chebyshev system, then $(g,f)$ is a negative one.

The most standard positive Chebyshev system on $\R$ is given by $f(x)=1$ and
$g(x)=x$. More generally, if $f,g:I\to\R$ are continuous functions with
$g\in\CP(I)$, $f/g\in\CM(I)$, then $(f,g)$ is a Chebyshev system. Indeed, we have
\Eq{D}{
  \D_{f,g}(x,y):=\DET{f(x)&f(y)\\g(x)&g(y)}
     =g(x)g(y)\left(\frac{f(x)}{g(x)}-\frac{f(y)}{g(y)}\right)
  \qquad(x,y\in I).
}
From, here it is obvious that $\D_{f,g}(x,y)$ vanishes if and only if $x=y$. Moreover, if $f/g$ is
decreasing (resp.\ increasing), then, for $x<y$, we have that $\D_{f,g}(x,y)>0$ (resp.\ $\D_{f,g}(x,y)<0$),
i.e., $(f,g)$ is a positive (resp.\ negative) Chebyshev system. By symmetry, analogous properties can be
established if $f$ is positive and $g/f$ strictly monotone.

For the sake of convenience and brevity, now we make the following hypotheses. We say that $m:I^d\times T\to
I$ is a \emph{measurable family of $d$-variable means on $I$} if
\begin{enumerate}[(H1)]
 \item $I$ is a nonvoid open real interval,
 \item $(T,\A)$ is a measurable space, where $\A$ is the $\sigma$-algebra of measurable sets of $T$,
 \item for all $t\in T$, $m(\cdot,t)$ is a $d$-variable mean on $I$,
 \item for all $\pmb{x}\in I^d$, the function $m(\pmb{x},\cdot)$ is measurable over $T$.
\end{enumerate}
If, instead of (H2) and H(4), we have that
\begin{enumerate}
 \item[(H2+)] $T$ is a topological space and $\A$ equals the $\sigma$-algebra $\B(T)$
of the Borel sets of $T$,
 \item[(H4+)] for all $\pmb{x}\in I^d$, the function $m(\pmb{x},\cdot)$ is continuous over $T$,
\end{enumerate}
then $m:I^d\times T\to I$ will be called a \emph{continuous family of $d$-variable means on $I$}.

For a measurable family of $d$-variable means $m:I^d\times T\to I$, we introduce the notations:
\Eq{*}{
  \underline{m}(\pmb{x}):=\inf_{t\in T} m(\pmb{x},t)
  \qquad\mbox{and}\qquad
  \overline{m}(\pmb{x}):=\sup_{t\in T} m(\pmb{x},t) \qquad(\pmb{x}\in I^d).
}
Obviously, by property (H3), for all $\pmb{x}\in I^d$, we have that $\min(\pmb{x})\leq
\underline{m}(\pmb{x})\leq\overline{m}(\pmb{x})\leq\max(\pmb{x})$. Provided that $T$ is a compact and
connected topological space and $m:I^d\times T\to I$ is a continuous family of $d$-variable means on $I$, we
have that
\Eq{mm}{
  [\underline{m}(\pmb{x}),\overline{m}(\pmb{x})]=\{m(\pmb{x},t)\mid t\in T\} \qquad(\pmb{x}\in I^d).
}

For the construction of a mean in terms of a Chebyshev system, a measurable family of means, and a
probability measure, we need the following basic lemma.

\Lem{0}{Let $m:I^d\times T\to I$ be a measurable family of $d$-variable means, let $\mu$ be a probability
measure on $(T,\A)$ and let $(f,g)$ be a Chebyshev system on $I$. Then, for all $\pmb{x}\in I^d$, there exists
a unique element $y\in [\underline{m}(\pmb{x}),\overline{m}(\pmb{x})]$ such that
\Eq{Eq}{
  \int_T \D_{f,g}(m(\pmb{x},t),y)\d\mu(t)=0.
}
Furthermore, if $(f,g)$ is a positive Chebyshev system, then, for all $u\in I$,
\Eq{*}{
  \int_T \D_{f,g}(m(\pmb{x},t),u)\d\mu(t)\lege 0 \qquad\mbox{if and only if}\qquad u\lege y.
}
In addition, if $g$ is positive and $f/g$ is strictly monotone, then
\Eq{Sol}{
   y=\left(\frac{f}{g}\right)^{-1}\left(
            \frac{\int_T f\big(m(\pmb{x},t)\big)\d\mu(t)}
                 {\int_T g\big(m(\pmb{x},t)\big)\d\mu(t)}\right).
}}

\begin{proof} Without loosing the generality, we may assume that $(f,g)$ is a positive Chebyshev system
throughout this proof.

For fixed $\pmb{x}\in I^d$, consider now the following function
\Eq{*}{
  h(u):=\int_T \D_{f,g}(m(\pmb{x},t),u)\d\mu(t)
       =g(u)\int_T f(m(\pmb{x},t))\d\mu(t)-f(u)\int_T g(m(\pmb{x},t))\d\mu(t) \qquad(u\in I).
}
By the continuity of $f$ and $g$, we have that $h$ is continuous on $I$. If $\overline{m}(\pmb{x})<u$,
then, for all $t\in T$ we have that $m(\pmb{x},t)<u$, hence $\D_{f,g}(m(\pmb{x},t),u)>0$. This implies that
$h(u)$ is positive for all $u\in I$ with $\overline{m}(\pmb{x})<u$. Similarly, for all $u\in I$ with
$u<\underline{m}(\pmb{x})$, we have that $h(u)<0$. Therefore, by the intermediate value property of
continuous functions, $h$ must have a zero between $\underline{m}(\pmb{x})$ and $\overline{m}(\pmb{x})$.

To prove the uniqueness, assume that $y$ and $z$ are distinct zeros of $h$ between
$\underline{m}(\pmb{x})$ and $\overline{m}(\pmb{x})$. Then we have that
\Eq{*}{
  g(y)\int_T f(m(\pmb{x},t))\d\mu(t)-f(y)\int_T g(m(\pmb{x},t))\d\mu(t)&=0,\\
  g(z)\int_T f(m(\pmb{x},t))\d\mu(t)-f(z)\int_T g(m(\pmb{x},t))\d\mu(t)&=0.
}
This means that the two unknowns $\xi:=\int_T f(m(\pmb{x},t))\d\mu(t)$ and $\eta:=\int_T
g(m(\pmb{x},t))\d\mu(t)$ are solutions of the following system of linear equations:
\Eq{*}{
  g(y)\xi-f(y)\eta&=0,\\
  g(z)\xi-f(z)\eta&=0.
}
Because $y$ and $z$ are distinct, the determinant of this system is nonzero, hence $\xi=\eta=0$. In this
case, $h(u)=0$ for all $u\in I$, which contradicts the property that $h(u)>0$ for
$\overline{m}(\pmb{x})<u$. The contradiction obtained shows that $y=z$, which proves the uniqueness of
the solution $y$ of equation \eq{Eq}. The uniqueness also implies that $h(u)>0$ for $u>y$ and $h(u)<0$ for
$u<y$.

Finally, assume that $g$ is positive and $f/g$ is strictly monotone (then $(f,g)$ is a Chebyshev system).
In this case, the equation $h(y)=0$ can be rewritten as
\Eq{*}{
   \frac{f(y)}{g(y)}=\frac{\int_T f(m(\pmb{x},t))\d\mu(t)}{\int_T g(m(\pmb{x},t))\d\mu(t)}.
}
By applying the inverse function of $f/g$ to this equation side by side, we obtain that $y$ is of the form
\eq{Sol}.
\end{proof}

The above lemma allows us to define a $d$-variable mean $M_{f,g,m;\mu}:I^d\to I$. Given $\pmb{x}\in I^d$, let
$M_{f,g,m;\mu}(\pmb{x})$ denote the unique solution $y$ of equation \eq{Eq}. In the particular case when $g$
is positive and $f/g$ is strictly monotone, we have that
\Eq{fg}{
   M_{f,g,m;\mu}(\pmb{x})
      :=\left(\frac{f}{g}\right)^{-1}\left(
            \frac{\int_T f\big(m(\pmb{x},t)\big)\d\mu(t)}
                 {\int_T g\big(m(\pmb{x},t)\big)\d\mu(t)}\right)
   \qquad(\pmb{x}\in I^d).
}
This mean will be called a $d$-variable \textit{generalized Bajraktarevi\'c mean} in the sequel.  When 
$g=1$, then 
\Eq{*}{
 M_{f,1,m;\mu}(\pmb{x})=(f)^{-1}\left(\int_T f\big(m(\pmb{x},t)\big)\d\mu(t)\right)
   \qquad(\pmb{x}\in I^d)
}
which will be termed a $d$-variable \textit{generalized quasi-arithmetic mean}.
If $T=\{1,\dots,d\}$, $\mu=\frac{\delta_1+\dots+\delta_d}{d}$ (where $\delta_{t}$ denotes the Dirac measure
concentrated at $t$) and $m(\pmb{x},t)=x_t$, then
\Eq{*}{
   M_{f,g,m;\mu}(\pmb{x})=M_{f,g}(\pmb{x})
    :=\left(\frac{f}{g}\right)^{-1}\bigg(\frac{f(x_1)+\dots+f(x_d)}{g(x_1)+\dots+g(x_d)}\bigg)
    \qquad\big(\pmb{x}=(x_1\dots,x_d)\in I^d\big),
}
which was introduced and studied by Bajraktarevi\'c \cite{Baj58}, \cite{Baj69}. When $g=1$, then 
$M_{f,1,m;\mu}(\pmb{x})=A_{f}(\pmb{x})$, which is the $d$-variable quasi-arithmetic mean introduced in
\eq{A}.

To define $d$-variable \textit{generalized Gini means}, let $p,q\in\R$ and assume that $I\subseteq\R_+$.
By taking
\Eq{pq}{
  f(x)&=x^p,&\quad g(x)&=x^q& \qquad\mbox{if}\quad p\neq q,\\
  f(x)&=x^p\log(x),&\quad g(x)&=x^p& \qquad\mbox{if}\quad p=q,
}
we can define $G_{p,q,m;\mu}$ in the following manner:
\Eq{ab}{
   G_{p,q,m;\mu}(\pmb{x})
      :=\begin{cases}
          \left(\dfrac{\int_T \big(m(\pmb{x},t)\big)^p\d\mu(t)}
                 {\int_T \big(m(\pmb{x},t)\big)^q\d\mu(t)}\right)^{\frac1{p-q}} &\mbox{if }p\neq q,\\[5mm]
          \exp\left(\dfrac{\int_T \big(m(\pmb{x},t)\big)^p\log\big(m(\pmb{x},t)\big)\d\mu(t)}
                 {\int_T \big(m(\pmb{x},t)\big)^p\d\mu(t)}\right) &\mbox{if }p=q.
        \end{cases}
   \qquad(\pmb{x}\in I^d).
}
In the particular case when $T=\{1,\dots,d\}$, $\mu=\frac{\delta_1+\dots+\delta_d}{d}$ and
$m(\pmb{x},t)=x_t$, the above formula reduces to the so-called $d$-variable \emph{Gini mean} $G_{p,q}$
(cf.\ \cite{Gin38}):
\Eq{*}{
   G_{p,q}(\pmb{x})
      :=\begin{cases}
         \bigg(\dfrac{x_1^p+\dots+x_d^p}{x_1^q+\dots+x_d^q}\bigg)^{\frac1{p-q}}&\mbox{if }p\neq q,\\[5mm]
         \exp\bigg(\dfrac{x_1^p\log(x_1)+\dots+x_d^p\log(x_d)}{x_1^p+\dots+x_d^p}\bigg)&\mbox{if }p=q,
        \end{cases}
    \qquad\big(\pmb{x}=(x_1\dots,x_d)\in\R_+^d\big).
}
Obviously, $G_{p,0}=H_p$, i.e., H\"older means are particular Gini means.

In what follows, we describe further interesting particular cases of formula \eq{fg}.
If $T=\{0,1,-1\}$, $I=\R_+$, $\mu=\mu_{-1}\delta_{-1}+\mu_{0}\delta_{0}+\mu_{1}\delta_{1}$ (where
$\mu_{-1}, \mu_{0}, \mu_{1}\in\R_+$ with $\mu_{-1}+\mu_{0}+\mu_{1}=1$), and $m(\pmb{x},t)=H_t(\pmb{x})$
(where $H_t$ stands for the $t$-th H\"older mean), then
\Eq{*}{
   M_{f,g,m;\mu}(\pmb{x})
      =\left(\frac{f}{g}\right)^{-1}\left(
            \frac{\mu_{-1}f(H_{-1}(\pmb{x}))+\mu_{0}f(H_{0}(\pmb{x}))+\mu_1f(H_{1}(\pmb{x}))}
                 {\mu_{-1}g(H_{-1}(\pmb{x}))+\mu_{0}g(H_{0}(\pmb{x}))+\mu_1g(H_{1}(\pmb{x}))}\right)
   \qquad(\pmb{x}\in\R_+^d).
}

In the next example we use the notations introduced in \eq{S} and \eq{At}.
If $T=S_d$, $\lambda$ is the $(d-1)$-dimensional Lebesgue measure on $S_d$,
and $m(\pmb{x},\pmb{t})=A_{\varphi}(\pmb{x},\pmb{t})$, where $\varphi\in\CM(I)$, then
\Eq{*}{
  M_{f,g,m;\mu}(\pmb{x})=M_{f,g,A_\varphi;\lambda}(\pmb{x})
      =\left(\frac{f}{g}\right)^{-1}\left(
            \frac{\int_{S_d}f\big(A_{\varphi}(\pmb{x},\pmb{t})\big)\d\lambda(\pmb{t})}
            {\int_{S_d}g\big(A_{\varphi}(\pmb{x},\pmb{t})\big)\d\lambda(\pmb{t})}\right)
   \qquad(\pmb{x}=(x_1,\dots,x_d)\in I^d).
}

If $\mu$ is the Lebesgue measure on $[0,1]$ and $f,g:I\to\R$ are
continuously differentiable functions such that $g'>0$ and $f'/g'$ is strictly monotone, and
$m(\pmb{x},t)=tx_1+(1-t)x_2$, then, by the Fundamental Theorem of Calculus, one can easily see that
\Eq{*}{
  M_{f',g', m; \mu}(\pmb{x})=C_{f,g}(\pmb{x})
   =\begin{cases}
        \bigg(\dfrac{f'}{g'}\bigg)^{-1}
             \bigg(\dfrac{f(x_2)-f(x_1)}{g(x_2)-g(x_1)}\bigg)
            &\mbox{if }x_1\neq x_2\\
         x_1  &\mbox{if }x_1=x_2
        \end{cases}
   \qquad(\pmb{x}=(x_1,x_2)\in I^2),
}
which is called a \emph{Cauchy} or \emph{difference mean} in the literature. Their equality problem was 
solved by Losonczi \cite{Los00a}.

By taking $f$ and $g$ given in \eq{pq}, the mean so obtained is the so-called \textit{Stolarsky mean},
which was introduced in the papers \cite{Sto75} and \cite{LeaSho78}. Their comparison problem was solved by 
Leach and Sholander \cite{LeaSho83} on unbounded intervals and by P\'ales \cite{Pal88b}, \cite{Pal92a} and by
Czinder and P\'ales \cite{CziPal06} on bounded intervals.

The aim of this paper is to study the \textit{global comparison problem}
\Eq{2}{
  M_{f,g,m;\mu}(\pmb{x})\le M_{h,k,n;\nu}(\pmb{x})\qquad(\pmb{x}\in I^d)
}
and also its \textit{local} analogue.
In terms of the Chebyshev systems $(f,g)$ and $(h,k)$, the measurable families of $d$-variable means
$m:I^d\times T\to I$ and $n:I^d\times S\to I$, and the measures $\mu,\nu$, we give necessary conditions
(which, in general, are not sufficient) and also sufficient conditions (that are also necessary in a certain
sense) for \eq{2} to hold. Our main results generalize that of the paper by Losonczi and P\'ales 
\cite{LosPal08} and also many former results obtained in various particular cases of this problem, cf.\ 
\cite{CziPal04}, \cite{CziPal05}, \cite{DarLos70}, \cite{Los02d}, \cite{Los02c}, \cite{NeuPal03}, 
\cite{Pal88c}.

\section{Invariants with respect to equality of means}

In order to describe the regularity conditions related to the two generating functions $f,g$ of the mean 
$M_{f,g,m;\mu}$, we introduce some regularity classes. The class $\C_0(I)$ consists of all those pairs of 
continuous functions $f,g:I\to\R$ that form a Chebyshev system over $I$.

If $k\ge1$, then we say that the pair $(f,g)$ is in the class $\C_k(I)$
if $f,g$ are $k$-times continuously differentiable functions such that $(f,g)\in\C_0(I)$ and the
\textit{Wronski determinant}
\Eq{W}{
  \DET{f'(x)&f(x)\\g'(x)&g(x)}
     =\partial_1\D_{f,g}(x,x)
   \qquad(x\in I)
}
does not vanish on $I$. Provided that $g$ is positive, then we have that
\Eq{f/g}{
  \bigg(\frac{f}{g}\bigg)'(x)=\frac{\partial_1\D_{f,g}(x,x)}{g^2(x)}
}
hence condition $\partial_1\D_{f,g}(x,x)\neq0$ implies that $f/g$ is
strictly monotone, whence it follows that $(f,g)\in\C_0(I)$.
Obviously, $\C_0(I)\supseteq\C_1(I)\supseteq\C_2(I)\supseteq\cdots$.

It is easy to see that if $(f,g),(f^*,g^*)\in\C_0(I)$ and
\Eq{5}{
  f&=\alpha f^*+\beta g^*,\\
  g&=\gamma f^*+\delta g^*,
}
where the constants $\alpha, \beta, \gamma, \delta\in \R$ satisfy
$\alpha\delta-\beta\gamma\ne 0$, then
\Eq{4.5}{
  \D_{f,g}(x,y)=\DET{\alpha&\beta\\\gamma&\delta}\cdot\D_{f^*,g^*}(x,y)
  \qquad(x,y\in I).
}
This implies that the identity
\Eq{4}{
  M_{f,g,m;\mu}=M_{f^*,g^*,m;\mu}
}
also holds for any measurable family $m:I^d\times T\to I$ and probability measure $\mu$.

If \eq{5} holds for some constants $\alpha, \beta, \gamma, \delta\in \R$, then we say that the pairs $(f,g)$
and $(f^*,g^*)$ are \emph{equivalent}. It is obvious that any necessary and/or sufficient condition for \eq{2}
has to be invariant with respect to the equivalence of the generating functions.

The following result, which is based on \cite[Theorem 3]{BesPal03}, allows us to assume more regularity on
Chebyshev systems.

\Lem{BP}{Let $k\in\N\cup\{0\}$ and $(f,g)\in\C_k(I)$. Then there exist $\alpha, \beta, \gamma, \delta\in \R$
with $\alpha\delta-\beta\gamma\ne 0$ and $(f^*,g^*)\in\C_k(I)$ such that \eq{5} holds and $g^*$ is positive
and $f^*/g^*$ is strictly monotone. Furthermore, if $k\geq1$, then the derivative of $f^*/g^*$ does not
vanish on $I$.}

\begin{proof} In the case $k=0$, the statement is a direct consequence of the result \cite[Theorem
3]{BesPal03} obtained by Bessenyei and P\'ales.

Assume now that $k\geq1$. Then, $(f,g)\in\C_k(I)$ means that $f$ and $g$ are $k$-times continuously
differentiable and $\partial_1\D_{f,g}(x,x)$ does not vanish for $x\in I$. Using what we have established
in the case $k=0$, there exist constants $\alpha, \beta, \gamma, \delta\in \R$ with
$\alpha\delta-\beta\gamma\ne 0$ and $(f^*,g^*)\in\C_0(I)$ such that \eq{5} holds and $g^*$ is positive and
$f^*/g^*$ is strictly monotone. The condition
$\alpha\delta-\beta\gamma\ne 0$ and \eq{5} imply that there exist $a,b,c,d\in\R$ with $ad-bc\neq0$ such that
\Eq{*}{
  f^*&=a f+b g,\\
  g^*&=c f+d g.
}
Hence $f^*$ and $g^*$ are also $k$-times continuously differentiable. This immediately implies that
\Eq{*}{
  \partial_1\D_{f^*,g^*}(x,x)=\DET{a&b\\c&d}\cdot\partial_1\D_{f,g}(x,x)
  \qquad(x\in I),
}
which shows that $\partial_1\D_{f^*,g^*}(x,x)$ does not vanish for $x\in I$. Therefore, $(f^*,g^*)\in\C_k(I)$
holds, too.

Applying the identity \eq{f/g} for $f^*$ and $g^*$ instead of $f$ and $g$, we can also see that the derivative
of $f^*/g^*$ does not vanish on $I$.
\end{proof}

For the computation of the first- and second-order partial derivatives of the mean $M_{f,g,m;\mu}$ at the
diagonal of $I^d$, we will establish a result below. For brevity, we introduce the following notation:
If $\pmb{p}\in I^d$ and $\delta>0$ then let $B(\pmb{p},\delta)$ stand for the ball $\{\pmb{x}\in I^d\colon
|\pmb{x}-\pmb{p}|\leq\delta\}$. Furthermore, if $\mu$ is a probability measure on the measurable
space $(T,\A)$ and $q\geq1$, then the space of measurable functions $\varphi:T\to\R$ such that
$|\varphi|^q$ is $\mu$-integrable will be denoted by $L^q(T,\A,\mu)$ or shortly by $L^q$.

If $\varphi,\psi:T\to\R$ are measurable functions such that $\varphi\psi$ is
$\mu$-integrable (for instance, $\varphi,\psi\in L^2$), then we set
\Eq{*}{
  \langle\varphi,\psi\rangle_\mu:=\int_T\varphi(t)\psi(t)\d\mu(t).
}
More generally, if $\varphi,\psi:I^d\times T\to\R$, and for some $\pmb{x}\in I^d$, the map
$t\mapsto\varphi(\pmb{x},t)\psi(\pmb{x},t)$ is $\mu$-integrable, then we write
\Eq{*}{
  \langle\varphi,\psi\rangle_\mu(\pmb{x}):=\int_T\varphi(\pmb{x},t)\psi(\pmb{x},t)\d\mu(t).
}
Given a number $q\geq1$, a function $\varphi:I^d\times T\to\R$ is said to be of $L^q$-type at $\pmb{p}\in
I^d$, if $\varphi(\pmb{p},\cdot)$ is measurable, furthermore, there exist $\delta>0$ and a function $a\in L^q$
such that
\Eq{*}{
  |\varphi(\pmb{x},t)|\leq a(t) \qquad(t\in T,\,\pmb{x}\in B(\pmb{p},\delta)).
}

Let $\C_1(I^d\times T)$ denote the class of measurable families of $d$-variable means $m:I^d\times T\to I$
with the following two additional properties:
\begin{enumerate}
 \item[(H5)] For every $t\in T$, the function $m(\cdot,t)$ is continuously partially differentiable over
$I^d$ such that, for all $\pmb{p}\in I^d$, $i\in\{1,\dots,d\}$, the function $\partial_i m$ is of
$L^1$-type at $\pmb{p}$.
\end{enumerate}
Analogously, we define $\C_2(I^d\times T)$ to be the following subclass of $\C_1(I^d\times T)$:
\begin{enumerate}
 \item[(H6)] For every $t\in T$, the function $m(\cdot,t)$ is twice continuously partially differentiable over
$I^d$ such that, for all $\pmb{p}\in I^d$ and $i,j\in\{1,\dots,d\}$, the function $\partial_i m$ is of
$L^2$-type and $\partial_i \partial_j m$ is of $L^1$-type at $\pmb{p}$.
\end{enumerate}

\Lem{1-}{Let $k\in\{1,2\}$ and let $\varphi:I\to\R$ be a $k$-times continuously differentiable function and
$m\in \C_k(I^d\times T)$. Then the function $\Phi:I^d\to\R$ defined by
\Eq{1.0}{
 \Phi(\pmb{x}):=\int_T\varphi(m(\pmb{x},t))\d\mu(t)
}
is $k$-times continuously differentiable on $I^d$. Furthermore, for $i\in\{1,\dots,d\}$,
\Eq{1.1}{
  \partial_i \Phi(\pmb{p})=\int_T\varphi'(m(\pmb{p},t))\,\partial_im(\pmb{p},t)\d\mu(t) \qquad(\pmb{p}\in I^d)
}
and, for $i,j\in\{1,\dots,d\}$,
\Eq{1.2}{
  \partial_i \partial_j \Phi(\pmb{p})
   =\int_T\Big[\varphi''(m(\pmb{p},t))\,\partial_im(\pmb{p},t)\,\partial_jm(\pmb{p},t)
   +\varphi'(m(\pmb{p},t))\,\partial_i\partial_jm(\pmb{p},t)\Big]\d\mu(t) \qquad(\pmb{p}\in I^d)
}
provided that $k=2$.}

\begin{proof} If $m$ is measurable family of $d$-variable means, then, by the continuity of $\varphi$ and the
mean value property of $m$, it easily follows that $\Phi$ is well-defined for all $\pmb{x}\in I^d$. Due to
the continuity of $\varphi'$ and the assumption $m\in \C_1(I^d\times T)$, it easily follows that the integral
on the right hand side of \eq{1.1} is well-defined. Furthermore, if $\varphi''$ is continuous and $m\in
\C_2(I^d\times T)$, then also the right hand side of \eq{1.2} exists.

First we elaborate the case $k=1$. We need to show that, for every $\pmb{p}\in I^d$, the
function $\Phi$ is partially differentiable at $\pmb{p}$ with formula \eq{1.1} and that the partial
derivatives are continuous at $\pmb{p}$.

Before proceeding to the proof, we shall establish, for every $\pmb{p}\in I^d$, the following equality
\Eq{Lim}{
  \lim_{\delta\to0}
  \int\limits_T\sup_{\pmb{x}\in B(\pmb{p},\delta)}\big|\varphi'(m(\pmb{x},t))\,\partial_im(\pmb{x},t)
       -\varphi'(m(\pmb{p},t))\,\partial_im(\pmb{p},t)\big|\d\mu(t)=0.
}
First choose $\delta_0>0$ so that $\alpha:=\min\{p_1,\dots,p_d\}-\delta_0$ and
$\beta:=\max\{p_1,\dots,p_d\}+\delta_0$ be elements of $I$. Let $K$ be the supremum of
$|\varphi'|$ over the compact interval $[\alpha,\beta]$. The continuity of $\varphi'$ implies that $K$ is
finite. By the mean value property of $m$, for every $t\in T$ and $\pmb{x}\in B(\pmb{p},\delta_0)$, we have
that $\alpha\leq m(\pmb{x},t)\leq\beta$. Hence, for every $t\in T$ and $\pmb{x}\in B(\pmb{p},\delta_0)$, the
inequality
$|\varphi'(m(\pmb{x},t))|\leq K$ holds.

Using the assumption that $\partial_im$ is of $L^1$-type at $\pmb{p}$, we can find
$0<\delta_1\leq\delta_0$ and a function $a\in L^1$ such that
\Eq{*}{
  |\partial_im(\pmb{x},t)|\leq a(t) \qquad(t\in T,\,\pmb{x}\in B(\pmb{p},\delta_1)).
}

Let $\delta_n>0$ be an arbitrary sequence converging to $0$ with $\delta_n\leq\delta_1$ for all $n\in\N$.
By the Lagrange Mean Value Theorem, for every $\pmb{x}\in B(\pmb{p},\delta_n)$ and for every $t\in T$,
there exists $\lambda\in[0,1]$ such that
\Eq{*}{
  |m(\pmb{x},t)-m(\pmb{p},t)|
    \leq\sum_{i=1}^d|\partial_im(\lambda\pmb{x}+(1-\lambda)\pmb{p},t)||x_i-p_i|
    \leq d\delta_n a(t).
}
Using the continuity of $\varphi'$ at $m(\pmb{p},t)$, it immediately follows that the sequence of measurable
functions $\psi_n:T\to\R$ defined by
\Eq{*}{
  \psi_n(t):=\sup_{\pmb{x}\in B(\pmb{p},\delta_n)}|\varphi'(m(\pmb{x},t))-\varphi'(m(\pmb{p},t))|\leq 2K
}
converges to zero for every $t\in T$. By the continuity of the partial derivative $\partial_im(\cdot,t)$ at
$\pmb{p}$, we also have that the sequence of measurable functions $\chi_n:T\to\R$ defined by
\Eq{*}{
  \chi_n(t):=\sup_{\pmb{x}\in B(\pmb{p},\delta_n)}|\partial_im(\pmb{x},t)-\partial_im(\pmb{p},t)|\leq 2a(t)
}
converges to zero for every $t\in T$. Using the above estimations, we can now obtain that
\Eq{*}{
  &\sup_{\pmb{x}\in B(\pmb{p},\delta_n)}\big|\varphi'(m(\pmb{x},t))\,\partial_im(\pmb{x},t)
       -\varphi'(m(\pmb{p},t))\,\partial_im(\pmb{p},t)\big|\\
  &\leq \sup_{\pmb{x}\in B(\pmb{p},\delta_n)}
       |\varphi'(m(\pmb{x},t))||\partial_im(\pmb{x},t)-\partial_im(\pmb{p},t)|
       +\sup_{\pmb{x}\in B(\pmb{p},\delta_n)}
       |\varphi'(m(\pmb{x},t))-\varphi'(m(\pmb{p},t))||\partial_im(\pmb{p},t)| \\
  &\leq  K \chi_n(t)+\psi_n(t) a(t).
}
The expression on the right hand side of this inequality converges to zero for each $t\in T$, and these
functions are dominated by the integrable function $4Ka$. Hence, by Lebesgue's Dominated Convergence Theorem,
it follows that
\Eq{*}{
  \lim_{n\to\infty}\int\limits_T\sup_{\pmb{x}\in
B(\pmb{p},\delta_n)}\big|\varphi'(m(\pmb{x},t))\,\partial_im(\pmb{x},t)
       -\varphi'(m(\pmb{p},t))\,\partial_im(\pmb{p},t)\big|=0.
}
Because the sequence $(\delta_n)$ converging to $0$ was arbitrary, it follows that \eq{Lim} holds.

Let $\pmb{p}\in I^d$ be fixed and let $\pmb{e}_i$ denote the $i$th vector of the standard basis on $\R^n$.
For the proof that the $i$th partial derivative of $\Phi$ at $\pmb{p}$ is given by \eq{1.1}, consider the
following estimation for $s\in (I-p_i)$, $s\neq0$:
\Eq{Delta}{
  \Delta_i(s):
  &= \bigg|\frac{\Phi(\pmb{p}+s\pmb{e}_i)-\Phi(\pmb{p})}{s}
    -\int_T\varphi'(m(\pmb{p},t))\,\partial_im(\pmb{p},t)\d\mu(t)\bigg| \\
&=\bigg|\frac1s\bigg(\int_T\varphi(m(\pmb{p}+s\pmb{e}_i,t))\d\mu(t)-\int_T\varphi(m(\pmb{p},t))\d\mu(t)\bigg)
  -\int_T\varphi'(m(\pmb{p},t))\,\partial_im(\pmb{p},t)\d\mu(t)\bigg| \\
  &\leq \frac1{|s|}\int_T\big|\varphi(m(\pmb{p}+s\pmb{e}_i,t))-\varphi(m(\pmb{p},t))
  -s\varphi'(m(\pmb{p},t))\,\partial_im(\pmb{p},t)\big|\d\mu(t).
}
Applying the Lagrange Mean Value Theorem for the function
\Eq{*}{
  \kappa_t(s):=\varphi(m(\pmb{p}+s\pmb{e}_i,t))-s\varphi'(m(\pmb{p},t))\,\partial_im(\pmb{p},t),
}
for every $t\in T$, we can find an element $\sigma_t$ between $0$ and $s$ such that
$\kappa_t(s)-\kappa_t(0)=\kappa_t'(\sigma_t)s$, that is
\Eq{*}{
  \varphi(m(\pmb{p}+s\pmb{e}_i,t))&-\varphi(m(\pmb{p},t))-s\varphi'(m(\pmb{p},t))\,\partial_im(\pmb{p},t) \\
  &=\big(\varphi'(m(\pmb{p}+\sigma_t\pmb{e}_i,t))\,\partial_im(\pmb{p}+\sigma_t\pmb{e}_i,t)
       -\varphi'(m(\pmb{p},t))\,\partial_im(\pmb{p},t)\big)s.
}
Using this formula, inequality \eq{Delta} and the equality \eq{Lim}, for $s\in I-p_i$, it follows that
\Eq{*}{
  \limsup_{s\to0} \Delta_i(s)
  &\leq \limsup_{s\to0} \int\limits_T
    \big|\varphi'(m(\pmb{p}+\sigma_t\pmb{e}_i,t))\,\partial_im(\pmb{p}+\sigma_t\pmb{e}_i,t)
       -\varphi'(m(\pmb{p},t))\,\partial_im(\pmb{p},t)\big|\d\mu(t) \\
   &\leq \lim_{s\to0} \int\limits_T \sup_{\pmb{x}\in B(\pmb{p},|s|)}
   \big|\varphi'(m(\pmb{x},t))\,\partial_im(\pmb{x},t)
       -\varphi'(m(\pmb{p},t))\,\partial_im(\pmb{p},t)\big|\d\mu(t)=0.
}
Thus, we have proved that $\Delta_i(s)$ tends to zero as $s\to0$. This completes the proof of the partial
differentiability of $\Phi$ with respect to the $i$th variable at $\pmb{p}$ and also the validity of formula
\eq{1.1}.

Finally, we show that the function $\partial_i\Phi$ is continuous at every $\pmb{p}\in I^d$.
Let $(\pmb{x}_n)$ be an arbitrary sequence in $B(\pmb{p},\delta_0)$ converging to $\pmb{p}$ and denote
$\delta_n:=|\pmb{x}_n-\pmb{p}|$. Then $(\delta_n)$ is a null sequence and we have that
\Eq{*}{
  |\partial_i\Phi(\pmb{x}_n)-\partial_i\Phi(\pmb{p})|
  &\leq \int\limits_T\big|\varphi'(m(\pmb{x}_n,t))\,\partial_im(\pmb{x}_n,t)
        -\varphi'(m(\pmb{p},t))\,\partial_im(\pmb{p},t)\big|\d\mu(t)\\
  &\leq \int\limits_T\sup_{\pmb{u}\in B(\pmb{p},\delta_n)}\big|\varphi'(m(\pmb{u},t))\,\partial_im(\pmb{u},t)
        -\varphi'(m(\pmb{p},t))\,\partial_im(\pmb{p},t)\big|\d\mu(t).
}
Due to the equality \eq{Lim}, the right hand side in the above inequality tends to zero as $n\to\infty$,
whence it follows that $(\partial_i\Phi(\pmb{x}_n))$ converges to $\partial_i\Phi(\pmb{p})$, which proves the 
continuity of $\partial_i\Phi$ at $\pmb{p}$.

Analogously, using a similar argument as in the proof of \eq{Lim}, for the case $k=2$, the reader can show
that the following two equalities hold:
\Eq{Lim1}{
  \lim_{\delta\to0}
  \int\limits_T\sup_{\pmb{x}\in B(\pmb{p},\delta)}\big|\varphi''(m(\pmb{x},t))\,\partial_im(\pmb{x},t)\,\partial_jm(\pmb{x},t)
       -\varphi''(m(\pmb{p},t))\,\partial_im(\pmb{p},t)\,\partial_jm(\pmb{p},t)\big|\d\mu(t)=0,
}
\Eq{Lim2}{
  \lim_{\delta\to0}
  \int\limits_T\sup_{\pmb{x}\in B(\pmb{p},\delta)}\big|\varphi'(m(\pmb{x},t))\,\partial_i\partial_jm(\pmb{x},t)
  -\varphi'(m(\pmb{p},t))\,\partial_i\partial_jm(\pmb{p},t)\big|\d\mu(t)=0.
}
Let $\pmb{p}\in I^d$ be fixed. To prove equality \eq{1.2} which establishes the formula for the $j$th partial
derivative of $\partial_i\Phi$ at $\pmb{p}$, consider the following estimation for $r\in (I-p_j)$, $r\neq0$:
\Eq{Omega}{
  \Delta_{ij}(r)&\\
  :=&\bigg|\frac{\partial_i\Phi(\pmb{p}+r\pmb{e}_j)-\partial_i\Phi(\pmb{p})}{r}
  -\int_T\varphi''(m(\pmb{p},t))\,\partial_im(\pmb{p},t)\,\partial_jm(\pmb{p},t)
  +\varphi'(m(\pmb{p},t))\,\partial_i\partial_jm(\pmb{p},t)\d\mu(t)\bigg|\\
  =&\bigg|\frac1r\bigg(\int_T\varphi'(m(\pmb{p}+r\pmb{e}_j,t))\,\partial_im(\pmb{p}+r\pmb{e}_j,t)\d\mu(t)-
  \int_T\varphi'(m(\pmb{p},t))\,\partial_im(\pmb{p},t)\d\mu(t)\bigg)\\
  &-\int_T\varphi''(m(\pmb{p},t))\,\partial_im(\pmb{p},t)\,\partial_jm(\pmb{p},t)
  +\varphi'(m(\pmb{p},t))\,\partial_i\partial_jm(\pmb{p},t)\d\mu(t)\bigg|\\
  \leq& \frac1{|r|}\int_T\big|\varphi'(m(\pmb{p}+r\pmb{e}_j,t))\,\partial_im(\pmb{p}+r\pmb{e}_j,t)-
  \varphi'(m(\pmb{p},t))\,\partial_im(\pmb{p},t)\\
  &-r\varphi''(m(\pmb{p},t))\,\partial_im(\pmb{p},t)\,\partial_jm(\pmb{p},t)
  -r\varphi'(m(\pmb{p},t))\,\partial_i\partial_jm(\pmb{p},t)\big|\d\mu(t).
}
Applying, for every $t\in T$, the Lagrange Mean Value Theorem for the function
\Eq{*}{
  \theta_t(r)\!:=\varphi'(m(\pmb{p}+r\pmb{e}_j,t))\partial_im(\pmb{p}+r\pmb{e}_j,t)
  -r\big(\varphi''(m(\pmb{p},t))\partial_im(\pmb{p},t)\partial_jm(\pmb{p},t)
  +\varphi'(m(\pmb{p},t))\partial_i\partial_jm(\pmb{p},t)\big),
}
we can find an element $\rho_t$ between $0$ and $r$ such that
\Eq{theta}{
\theta_t(r)-\theta_t(0)=\theta_t'(\rho_t)r.
}
Now, by using equality \eq{theta}, inequality \eq{Omega} and equalities \eq{Lim1} and \eq{Lim2},
respectively, with an analogous argument that we applied in the case $k=1$, we get that $\Delta_{ij}(r)$
tends to zero as $r\to0$, proving the partial differentiability of $\partial_i\Phi$ at
$\pmb{p}$ with respect to the $j$th variable and formula \eq{1.2}. On the other hand, again by a similar
train of thoughts, it easily follows from  \eq{Lim1} and \eq{Lim2} that the function
$\partial_i\partial_j\Phi$ is continuous on $I^d$. This completes the proof of the lemma.
\end{proof}

\Thm{1+}{Let $(f,g)\in\C_1(I)$, let $m\in\C_1(I^d\times T)$ be a measurable family of means, and let $\mu$ be
a probability measure on the measurable space $(T,\A)$.
Then $M_{f,g,m;\mu}$ is continuously differentiable on $I^d$ and, for all $i\in\{1,\dots,d\}$ and $x\in I$,
\Eq{PD1}{
  \partial_iM_{f,g,m;\mu}(x,\dots,x)
  =\langle\partial_im,1\rangle_\mu(x,\dots,x).
}
If, in addition, $(f,g)\in\C_2(I)$, let $m\in\C_2(I^d\times T)$, then $M_{f,g,m;\mu}$ is twice continuously
differentiable on $I^d$ and, for all $i,j\in\{1,\dots,d\}$ and $x\in I$,
\Eq{PD2}{
  \partial_i\partial_j &M_{f,g,m;\mu}(x,\dots,x)\\
  &=\big(\langle \partial_im,\partial_j m\rangle_\mu
    -\langle\partial_im,1\rangle_\mu\langle\partial_jm,1\rangle_\mu\big)(x,\dots,x)
    \frac{\partial_1^2\D_{f,g}(x,x)}{\partial_1\D_{f,g}(x,x)}
    +\langle \partial_i\partial_j m,1\rangle_\mu(x,\dots,x).
}}

\begin{proof} Let $k\in\{1,2\}$ and assume that $(f,g)\in\C_k(I)$, $m\in\C_k(I^d\times T)$. In view of
\lem{BP}, we may assume that $g$ is positive, $f/g$ is strictly monotone with a non-vanishing first-order
derivative. Then $f$, $g$ and the inverse of $f/g$ are $k$-times continuously differentiable and, by \lem{1-},
we also have that the mappings
\Eq{*}{
  \pmb{x}\mapsto \int_T f\big(m(\pmb{x},t)\big)\d\mu(t)
  \qquad\mbox{and}\qquad
  \pmb{x}\mapsto \int_T g\big(m(\pmb{x},t)\big)\d\mu(t)
}
are $k$-times continuously differentiable on $I^d$. On the other hand, we now also
have formula \eq{fg} for the $d$-variable mean $M_{f,g,m;\mu}$. Thus, using the standard calculus rules, it
follows that $M_{f,g,m;\mu}$ is $k$-times continuously differentiable on $I^d$

To prove the first formula stated in \eq{PD1}, let us consider the case $k=1$.
Differentiating the identity \eq{Eq} with respect to the variable $x_i$ once, we get
\Eq{*}{
   \int_T \big[\partial_1\D_{f,g}(m(\pmb{x},t),M_{f,g,m;\mu}(\pmb{x}))\partial_im(\pmb{x},t)
+\partial_2\D_{f,g}(m(\pmb{x},t),M_{f,g,m;\mu}(\pmb{x}))\partial_iM_{f,g,m;\mu}(\pmb{x})\big]\d\mu(t)=0.
}
Now taking $x\in I$ and substituting $\pmb{x}=(x,\dots,x)$, the above equation simplifies to
\Eq{*}{
   \partial_1\D_{f,g}(x,x)\int_T \partial_im(x,\dots,x,t)\d\mu(t)
+\partial_2\D_{f,g}(x,x)\partial_iM_{f,g,m;\mu}(x,\dots,x)=0.
}
Observe that $\partial_1\D_{f,g}(x,x)=-\partial_2\D_{f,g}(x,x)\neq0$ for all $x\in I$, hence the former
equation yields the desired equality \eq{PD1}.

Now consider the case $k=2$. Differentiating the identity that we obtained in the first
lines of the proof with respect to the variable $x_j$, we obtain
\Eq{*}{
   &\int_T \big[\partial_1^2\D_{f,g}(m(\pmb{x},t),M_{f,g,m;\mu}(\pmb{x}))
             \partial_im(\pmb{x},t)\partial_jm(\pmb{x},t) \\
 &+\partial_1\partial_2\D_{f,g}(m(\pmb{x},t),M_{f,g,m;\mu}(\pmb{x}))
             \big(\partial_jM_{f,g,m;\mu}(\pmb{x})\partial_im(\pmb{x},t)
             +\partial_iM_{f,g,m;\mu}(\pmb{x})\partial_jm(\pmb{x},t)\big)\\
&+\partial_2^2\D_{f,g}(m(\pmb{x},t),M_{f,g,m;\mu}(\pmb{x}))
             \partial_jM_{f,g ,m;\mu}(\pmb{x})\partial_iM_{f,g ,m;\mu}(\pmb{x})\\
&+\partial_1\D_{f,g}(m(\pmb{x},t),M_{f,g,m;\mu}(\pmb{x}))
             \partial_j\partial_im(\pmb{x},t)
 +\partial_2\D_{f,g}(m(\pmb{x},t),M_{f,g,m;\mu}(\pmb{x}))
             \partial_j\partial_iM_{f,g,m;\mu}(\pmb{x})\big]\d\mu(t)=0 ,
}
respectively. Using the identities
$\partial_2\D_{f,g}(x,x)=-\partial_1\D_{f,g}(x,x)$,
$\partial_2^2\D_{f,g}(x,x)=-\partial_1^2\D_{f,g}(x,x)$,
and $\partial_1\partial_2\D_{f,g}(x,x)=0$ (that are consequences of the
asymmetry property $\D_{f,g}(x,y)=-\D_{f,g}(y,x)$), and substituting $\pmb{x}=(x,\dots,x)$, we get that
\Eq{*}{
   & \partial_1^2\D_{f,g}(x,x)\int_T \partial_im(x,\dots,x,t)\partial_jm(x,\dots,x,t)\d\mu(t)
   +\partial_1\D_{f,g}(x,x)\int_T\partial_j\partial_im(x,\dots,x,t)\d\mu(t)\\
&-\partial_1^2\D_{f,g}(x,x)
             \partial_jM_{f,g ,m;\mu}(x,\dots,x)\partial_iM_{f,g ,m;\mu}(x,\dots,x)
 -\partial_1\D_{f,g}(x,x)\partial_j\partial_iM_{f,g,m;\mu}(x,\dots,x)=0.
}
Dividing both sides of this equation by $\partial_1\D_{f,g}(x,x)\neq0$ and using \eq{PD1}, for the
second-order partial derivative $\partial_j\partial_iM_{f,g,m;\mu}(x,\dots,x)$, we obtain the formula stated
in \eq{PD2}.
\end{proof}

One of the most important particular case of the above theorem is when the $d$-variable family of means is 
a family of weighted $d$-variable arithmetic means.

\Cor{1+}{Let $(f,g)\in\C_1(I)$, let $\mu$ be a probability measure on the measurable space $(T,\A)$ and let 
$m\in\C_1(I^d\times T)$ be a measurable family of $d$-variable means given by
\Eq{*}{
  m(\pmb{x},t):=\lambda_1(t)x_1+\cdots+\lambda_d(t)x_d \qquad(\pmb{x}=(x_1,\dots,x_d)\in I^d,\,t\in T),
}
where $\lambda_1,\dots,\lambda_d:T\to[0,1]$ are measurable functions with $\lambda_1+\cdots+\lambda_d=1$.
Then $M_{f,g,m;\mu}$ is continuously differentiable on $I^d$ and, for all $i\in\{1,\dots,d\}$ and $x\in I$,
\Eq{PD1+}{
  \partial_iM_{f,g,m;\mu}(x,\dots,x)=\langle\lambda_i,1\rangle_\mu.
}
If, in addition, $(f,g)\in\C_2(I)$, let $m\in\C_2(I^d\times T)$, then $M_{f,g,m;\mu}$ is twice continuously
differentiable on $I^d$ and, for all $i,j\in\{1,\dots,d\}$ and $x\in I$,
\Eq{PD2+}{
  \partial_i\partial_j &M_{f,g,m;\mu}(x,\dots,x)
  =\big(\langle\lambda_i,\lambda_j\rangle_\mu
    -\langle\lambda_i,1\rangle_\mu\langle\lambda_j,1\rangle_\mu\big)
    \frac{\partial_1^2\D_{f,g}(x,x)}{\partial_1\D_{f,g}(x,x)}.
}}

\begin{proof} Observe that we have $\partial_im (\pmb{x},t)=\lambda_i(t)$ and 
$\partial_i\partial_jm(\pmb{x},t)=0$. By the boundedness of the measurable function, it follows that 
$\partial_im$ is $L^1$- and $L^2$-type, and $\partial_i\partial_jm$ is $L^1$-type at every point of $I^d$.
Therefore, \thm{1+} applies, and formulas \eq{PD1} and \eq{PD2} reduce to \eq{PD1+} and \eq{PD2+}, 
respectively.
\end{proof}

The particular case when $T=[0,1]$, $d=2$, $\lambda_1(t)=t$ and $\lambda_2=1-t$ was considered by Losonczi 
and P\'ales in the paper \cite{LosPal08}, where also the related local and global comparison problems were 
investigated. The above \thm{1+} and \cor{1+} generalize the result of \cite[Lemma 4]{LosPal08}.

The following lemma, which is an extension of \cite[Lemma 3]{LosPal08}, will play an important role in 
establishing the necessary conditions for the (global) comparison of means. We recall that a sequence 
$(\nu_k)$ of probability measures on $T$ is said to \emph{converge weakly to a measure} $\nu$ if, for all 
bounded Borel measurable functions $\varphi:T\to\R$, we have
\Eq{*}{
  \lim_{k\to\infty}\int_T \varphi(t)\d\nu_k(t)=\int_T \varphi(t)\d\nu(t).
}

\Lem{2}{Let $(f,g)\in\C_1(I)$ and let $(\nu_k)$ be a sequence of probability measures on $T$ weakly
converging to a measure $\nu$, let $(\gamma_k)$ be a null sequence of positive numbers in $[0,1]$ and let
$t_0\in T$. Set $\mu_k:=(1-\gamma_k)\delta_{t_0}+\gamma_k\nu_k$ for $k\in\N$. Then
\Eq{G}{
   \lim_{k\to\infty}\frac1\gamma_k&\big[M_{f,g,m;\mu_k}(\pmb{x})-m(\pmb{x},t_0)\big]
    = \frac{\int_T
\D_{f,g}(m(\pmb{x},t),m(\pmb{x},t_0))\d\nu(t)}{\partial_1\D_{f,g}(m(\pmb{x},t_0),m(\pmb{x},t_0))}
\qquad(\pmb{x}\in I^{d}).
}}

\begin{proof} Let $\pmb{x}\in I^d$ be a fixed vector. By the assumptions of the lemma, we have that $\mu_k$
converges to $\delta_{t_0}$ weakly. More generally, for an arbitrary bounded sequence of Borel measurable
functions $\varphi_{k}:T\to\R$, which converges uniformly to $\varphi_0:T\to\R$ as $k\to\infty$, we get
\Eq{F1}{
  \int_T \varphi_{k} \d\mu_{k}(t) & =\int_T \varphi_{k} d((1-\gamma_{k})\delta_{t_0}+\gamma_{k}\nu_{k})(t) \\
   & = (1-\gamma_{k})\int_T \varphi_{k} d\delta_{t_0}(t) + \gamma_{k}\int_T \varphi_0 \d\nu_{k}(t)
     +\gamma_{k}\int_T (\varphi_{k}-\varphi_0) \d\nu_{k}(t)\\
   &\longrightarrow  \int_T \varphi_0 d\delta_{t_0}(t)=\varphi_0(t_0).
}
First, we are going to show that the sequence $u_{k}:=M_{f,g,m;\mu_k}(\pmb{x})$ converges to $m(\pmb{x},t_0)$.
We have that $\min(\pmb{x})\leq u_k\leq\max(\pmb{x})$ for all $k\in\N$. Hence it is sufficient to prove that 
every convergent subsequence of $(u_k)$ converges to the same limit point. To show this, let $(u_{k_j})$ be 
any convergent subsequence of $(u_k)$ such that $u_{k_j}\to u_{0}$ as $j\to\infty$. Then, the sequence of 
Borel measurable functions $\varphi_j(t):=D_{f,g}(m(\pmb{x},t), u_{k_j})$ tends uniformly to the limit 
function $\varphi_0(t):=D_{f,g}(m(\pmb{x},t), u_{0})$. Thus, in view of formula \eq{F1}, we get
\Eq{*}{
  \lim_{j\to\infty}\int_T D_{f, g}(m(\pmb{x},t),u_{k_j})\d\mu_{k_j}(t)=D_{f, g}(m(\pmb{x},t_0), u_{0}).
}
On the other hand, for all $j$, we have that
\Eq{*}{
  \int_T D_{f, g}(m(\pmb{x},t), u_{k_j})\d\mu_{k_j}(t)=0,
}
which implies that $D_{f, g}(m(\pmb{x},t_0), u_{0})$ is zero, i.e., $u_0=m(\pmb{x},t_0)$.
Hence $u_{k}\to m(\pmb{x},t_0)$ as $k\to\infty$.

Moreover, as $k\to\infty$, we similarly obtain
\Eq{LF}{
\frac{1}{\gamma_{k}}&\int_T D_{f, g}(m(\pmb{x},t), m(\pmb{x},t_0))\d\mu_{k}(t)\\
  &=\frac{1}{\gamma_{k}}\int_T D_{f, g}(m(\pmb{x},t),
m(\pmb{x},t_0))\d((1-\gamma_{k})\delta_{t_0}+\gamma_{k}\nu_{k})(t)\\
  &=\int_T D_{f, g}(m(\pmb{x},t), m(\pmb{x},t_0))\d\nu_{k}(t)\longrightarrow\int_T D_{f,
g}(m(\pmb{x},t),
m(\pmb{x},t_0))\d\nu(t).
}

Taking $\Phi_k(u):=\int_T D_{f, g}(m(\pmb{x},t), u)\d\mu_{k}(t)$ and applying the Lagrange mean value theorem
for the differentiable function $\Phi_k$, for every $k\in\N$, we can find a number $\eta_k$ between $u_{k}$
and $m(\pmb{x},t_0)$ such that
\Eq{LMVT3}{
 \Phi_k(m(\pmb{x},t_0))-\Phi_k(u_{k})=\Phi_k'(\eta_{k})(m(\pmb{x},t_0)-u_{k})
}
Since $\Phi_k(u_{k})=0$, therefore it follows that
\Eq{LMVT2}{
  \frac{1}{\gamma_{k}}\Phi_k(m(\pmb{x},t_0))
  = \Phi_k'(\eta_{k})\cdot\frac{1}{\gamma_{k}}(m(\pmb{x},t_0)-u_{k}).
}
Thus,
\Eq{LMVT}{
\frac{1}{\gamma_{k}}&\int_T D_{f, g}(m(\pmb{x},t), m(\pmb{x},t_0))\d\mu_{k}(t)\\
  &\hspace{1cm}
  =\int_T\DET{f(m(\pmb{x},t))&f'(\eta_{k})\\g(m(\pmb{x},t))&g'(\eta_{k})}\d\mu_{k}(t)
  \cdot \frac{1}{\gamma_{k}}\big[m(\pmb{x},t_0)-M_{f,g,m;\mu_k}(\pmb{x})\big].
}
Then, obviously, $\eta_{k}$ converges to $m(\pmb{x},t_0)$. By taking the limit of both sides of \eq{LMVT} as
$k\to\infty$ and using \eq{LF} and \eq{F1}, we get
\Eq{LMVT1}{
\int_T
     &\D_{f,g}(m(\pmb{x},t),m(\pmb{x},t_0))\d\nu(t)\\[-3mm]
      &\hspace{1cm}=\partial_1\D_{f,g}(m(\pmb{x},t_0),m(\pmb{x},t_0))\cdot
\lim_{k\to\infty}\frac1\gamma_k\big[M_{f,g,m;\mu_k}(\pmb{x})-m(\pmb{x},t_0)\big]
}
By dividing both sides by $\partial_1\D_{f,g}(m(\pmb{x},t_0),m(\pmb{x},t_0))$, we get
\Eq{LMVT4}{
\lim_{k\to\infty}\frac1\gamma_k & [M_{f,g,m;\mu_k}(\pmb{x})-m(\pmb{x},t_0)]=\frac{\int_T
      \D_{f,g}(m(\pmb{x},t),m(\pmb{x},t_0))\d\nu(t)}{\partial_1\D_{f,g}(m(\pmb{x},t_0),m(\pmb{x},t_0))}.}
This completes the proof of the lemma.
\end{proof}

\section{Necessary conditions, sufficient conditions for local comparison of means}

Our first result offers a necessary as well as a sufficient condition for the local comparison of means.
Given two $d$-variable means $M,N:I^d\to I$, we say that \emph{$M$ is locally smaller than $N$ at
$x_0\in I$} if there exists a neighborhood $U\subseteq I$ of $x_0$ such that
\Eq{1C}{
  M(\pmb{x})\leq N(\pmb{x})
}
holds for all $\pmb{x}\in U^d$. The case $d=1$ being trivial, we always assume that $d\geq2$ holds in the
subsequent considerations.

\Thm{1}{Let $M,N:I^d\to I$ be $d$-variable means such that $M$ is locally smaller than $N$ at a point $x_0\in
I$. Assume that $M$ and $N$ are partially differentiable at the diagonal point $(x_0,\dots,x_0)\in I^d$.
Then, for $x=x_0$ and for all $i\in\{1,\dots,d\}$,
\Eq{1A}{
  \partial_i M(x,\dots,x)=\partial_i N(x,\dots,x).
}
If, in addition, $M$ and $N$ are twice differentiable at $(x_0,\dots,x_0)\in I^d$, then the
symmetric $(d-1)\times (d-1)$-matrix
\Eq{1B}{
  \big(\partial_i\partial_j N(x_0,\dots,x_0)-\partial_i\partial_j M(x_0,\dots,x_0)\big)_{i,j=1}^{d-1}
}
is positive semidefinite. \\\indent On the other hand, if, for some $x_0\in I$, the equality \eq{1A} holds for
all $i\in\{1,\dots,d\}$ and for all $x$ in a neighborhood of $x_0$, furthermore, $M$ and $N$ are twice
continuously differentiable at $(x_0,\dots,x_0)$ and the symmetric $(d-1)\times (d-1)$-matrix given by \eq{1B}
is positive definite, then $M$ is locally smaller than $N$ at $x_0$.}

\begin{proof} Assume that $M$ is locally smaller than $N$ at $x_0\in I$, i.e., \eq{1C} holds for all
$\pmb{x}\in U^d$ in a neighborhood $U\subseteq I$ of $x_0$. Assume that $M$ and $N$ are partially
differentiable at the diagonal point $(x_0,\dots,x_0)\in I^d$. Define the function $D:U^d\to\R$ by
\Eq{*}{
  D(\pmb{x})= N(\pmb{x})-M(\pmb{x}) \qquad (\pmb{x}\in U^d).
}
Then $D$ is nonnegative by inequality \eq{1C} and attains its minimum (which equals zero) at
$\pmb{x}=(x_0,\dots,x_0)$. Therefore $\partial_iD(x_0,\dots,x_0)=0$ for all $i\in\{1,\dots,d\}$, which yields
\eq{1A}.

If, in addition, $M$ and $N$ are twice differentiable at $(x_0,\dots,x_0)\in I^d$. Then
$D''(x_0,\dots,x_0)=\big(\partial_i\partial_j D(x_0,\dots,x_0)\big)_{i,j=1}^d$ is a positive semidefinite
symmetric $d\times d$-matrix. By the well-known necessary and sufficient conditions of positive
semidefiniteness (cf.\ \cite{Bha07}), this implies that the symmetric $(d-1)\times(d-1)$-matrix
$\big(\partial_i\partial_j D(x_0,\dots,x_0)\big)_{i,j=1}^{d-1}$ is also positive semidefinite.

Now let $x_0\in I$ and assume that there exists a neighborhood $U\subseteq I$ of $x_0$ such that the means $M$
and $N$ are twice differentiable on $U^d$, their second-order partial derivatives are continuous at
$\pmb{x}_0=(x_0,\dots,x_0)$, the equality \eq{1A} holds for all $x\in U$ and for
all $i\in\{1,\dots,d\}$, and the symmetric $(d-1)\times (d-1)$-matrix-valued function
$A(\pmb{x}):=\big(\partial_i\partial_j D(\pmb{x})\big)_{i,j=1}^{d-1}$ is positive definite at
$\pmb{x}_0$.

By Sylvester's criterion, $A(\pmb{x}_0)$ is positive definite if and only if all of its leading
principal minors are positive. By the continuity of the second-order partial derivatives, $A$ is continuous
at $\pmb{x}_0$, hence its leading principal minors are also continuous at $\pmb{x}_0$. Therefore, there is a
neighborhood $V\subseteq I^d$ of $\pmb{x}_0$ where these leading principal minors are positive and hence, at
the points of $V$, $A$ is also positive definite. By shrinking the neighborhood $U$ of $x_0$ if necessary, we
may assume that $U$ is an interval and $U^d\subseteq V$. Hence $A(\pmb{x})$ is positive definite for all
$\pmb{x}\in U^d$.

In order to show that the inequality \eq{1C} holds for all $\pmb{x}\in U^d$, let $\pmb{x}=(x_1,\dots,x_d)\in
U^d$ be fixed and apply the Taylor Mean Value Theorem to the function
\Eq{*}{
  (u_1,\dots,u_{d-1})\mapsto D(u_1,\dots,u_{d-1},x_d) \qquad\big((u_1,\dots,u_{d-1})\in U^{d-1}\big)
}
at the base point $(x_d,\dots,x_d)\in U^{d-1}$. In view of this theorem, there exists $\theta\in[0,1]$, such
that
\Eq{1D}{
  D(x_1,\dots, & x_{d-1},x_d) \\
  =& D(x_d,\dots,x_d,x_d)+\sum_{i=1}^{d-1} \partial_iD(x_d,\dots,x_d,x_d)(x_i-x_d) \\
   & +\frac12\sum_{i=1}^{d-1}\sum_{j=1}^{d-1}\partial_i\partial_j
     D(\theta x_1+(1-\theta)x_d,\dots,\theta x_{d-1}+(1-\theta)x_d,x_d)(x_i-x_d)(x_j-x_d).
}
We have that $D(x_d,\dots,x_d,x_d)=M(x_d,\dots,x_d,x_d)-N(x_d,\dots,x_d,x_d)=x_d-x_d=0$, equation \eq{1A}
applied for $x=x_d$ implies that $\partial_iD(x_d,\dots,x_d,x_d)=0$ for all $i\in\{1,\dots,d-1\}$. Finally,
$A$ is positive definite at the point $(\theta x_1+(1-\theta)x_d,\dots,\theta x_{d-1}+(1-\theta)x_d,x_d)\in
U^d$, hence the last term on the right hand side of \eq{1D} is nonnegative. Thus \eq{1D} shows that
$D(x_1,\dots,x_{d-1},x_d)\geq0$, which implies that $M$ is smaller than $N$ on $U^d$.
\end{proof}

\Rem{1}{We note that, for the sufficiency part of the theorem, the standard 2nd-order sufficient condition
for the local minimum cannot be applied. The reason is that the matrix
\Eq{MN}{
  \big(\partial_i\partial_j N(\pmb{x}_0)-\partial_i\partial_j M(\pmb{x}_0)\big)_{i,j=1}^{d}
}
can never be positive definite. Indeed, if $M$ is locally smaller than $N$ at $x_0$, then $M$ is locally
smaller than $N$ at every $x$ in a neighborhood $U$ of $x_0$ and hence \eq{1A} holds for all $x\in U$ and
$i\in\{1,\dots,d\}$. Differentiating \eq{1A} with respect to $x$ at $x_0$, we obtain, for all
$i\in\{1,\dots,d\}$, that
\Eq{*}{
  \sum_{j=1}^d \partial_j\partial_i N(\pmb{x}_0)=\sum_{j=1}^d \partial_j\partial_i M(\pmb{x}_0).
}
This shows that the sum of the columns of the matrix in \eq{MN} is the zero vector. Therefore, the
determinant of this matrix is zero, showing that this matrix is not positive definite.}

\Cor{1}{Let $(f,g),(h,k)\in\C_1(I)$, let $m\in\C_1(I^d\times T)$ and $n\in\C_1(I^d\times S)$ be measurable
families of means, and let $\mu$ and $\nu$ be probability measures on the measurable spaces $(T,\A)$ and
$(S,\B)$, respectively. Suppose that $M_{f,g,m;\mu}$ is locally smaller than $M_{h,k,n;\nu}$ at $x_0\in I$.
Then, there exists a neighborhood $U\subseteq I$ of $x_0$ such that for $x\in U$ and for all
$i\in\{1,\dots,d\}$,
\Eq{6}{
  \langle\partial_im,1\rangle_\mu(x,\dots,x)=\langle\partial_in,1\rangle_\nu(x,\dots,x).
}
If, in addition, $(f,g),(h,k)\in\C_2(I)$, $m\in\C_2(I^d\times T)$, and $n\in\C_2(I^d\times S)$, then the
$(d-1)\times(d-1)$-matrix whose $(i,j)$th entry is given by
\Eq{7}{
  &\big(\langle \partial_in,\partial_j n\rangle_\nu
    -\langle\partial_in,1\rangle_\nu\langle\partial_jn,1\rangle_\nu\big)(\pmb{x}_0)
    \frac{\partial_1^2\D_{h,k}(x_0,x_0)}{\partial_1\D_{h,k}(x_0,x_0)}
    +\langle \partial_i\partial_j n,1\rangle_\nu(\pmb{x}_0)\\
    &-\big(\langle \partial_im,\partial_j m\rangle_\mu
    -\langle\partial_im,1\rangle_\mu\langle\partial_jm,1\rangle_\mu\big)(\pmb{x}_0)
    \frac{\partial_1^2\D_{f,g}(x_0,x_0)}{\partial_1\D_{f,g}(x_0,x_0)}
    -\langle \partial_i\partial_j m,1\rangle_\mu(\pmb{x}_0)
}
for $i,j\in\{1,\dots,d-1\}$ is positive semidefinite. \\\indent On the other hand, if $(f,g),(h,k)\in\C_2(I)$,
$m\in\C_2(I^d\times T)$, $n\in\C_2(I^d\times S)$, and \eq{6} holds for all $i\in\{1,\dots,d\}$ and for all
$x$ in a neighborhood of $x_0$ and the $(d-1)\times (d-1)$-matrix whose $(i,j)$th entry is given by \eq{7} is
positive definite, then $M_{f,g,m;\mu}$ is locally smaller than $M_{h,k,n;\nu}$ at $x_0\in I$. }

\begin{proof}
If $(f,g),(h,k)\in\C_k(I)$, $m\in\C_k(I^d\times T)$, and $n\in\C_k(I^d\times S)$, then \thm{1+} implies that
$M_{f,g,m;\mu}$ and $M_{h,k,n;\nu}$ are $k$-times continuously differentiable on $I^d$ in the cases
$k\in\{1,2\}$.

Assume that $M_{f,g,m;\mu}$ is locally smaller than $M_{h,k,n;\nu}$ at $x_0\in I$. Then, by \thm{1}, there
exists a neighborhood $U\subseteq I$ of $x_0$ such that for $x\in U$ and for all $i\in\{1,\dots,d\}$,
\Eq{*}{
  \partial_iM_{f,g,m;\mu}(x,\dots, x)=\partial_iM_{h,k,n;\nu}(x,\dots, x).
}
Applying formula \eq{PD1} of \thm{1+}, the necessity of condition \eq{6} follows.

In addition, if the second-order regularity assumptions are satisfied, then, by \thm{1}, the
$(d-1)\times(d-1)$-matrix whose $(i,j)$th entry is given by
\Eq{*}{
  \big(\partial_i\partial_j M_{h,k,n;\nu}(x_0,\dots,x_0)-\partial_i\partial_j M_{f,g,m;\mu}(x_0,\dots,x_0)\big)_{i,j=1}^{d-1}
}
for $i,j\in\{1,\dots,d-1\}$ is positive semidefinite. Now the application of formula \eq{PD2} of \thm{1+}
yields the necessity of condition \eq{7}.

Now, under the second-order regularity assumptions suppose that \eq{6} holds for all $i\in\{1,\dots,d\}$ and
for all $x$ in a neighborhood of $x_0$ and the $(d-1)\times (d-1)$-matrix whose $(i,j)$th entry is given by
\eq{7} is positive definite. Since $M_{f,g,m;\mu}$ and $M_{h,k,n;\nu}$ are twice continuously differentiable
and by \thm{1}, we have $M_{f,g,m;\mu}$ is locally smaller than $M_{h,k,n;\nu}$ at $x_0\in I$.
\end{proof}

In the special setting when $T=[0,1]$, $d=2$, $m$ is given by $m(\pmb{x},t):=tx_1+(1-t)x_2$, the above 
\cor{1} simplifies to the result of \cite[Theorem 5]{LosPal08}. Now we consider the particular case when the 
families of means $m$ and $n$ as well as the measures $\mu$ and $\nu$ coincide.

\Cor{2}{Let $(f,g),(h,k)\in\C_2(I)$, let $m\in\C_2(I^d\times T)$ be a measurable family of means, and let
$\mu$ be a probability measure on the measurable space $(T,\A)$. Let $x_0\in I$ and assume that there exists
$i\in\{1,\dots,d-1\}$ such that, the map $t\mapsto\partial_im(\pmb{x_0},t)$ is not $\mu$-almost everywhere
constant on $T$. If $M_{f,g,m;\mu}$ is locally smaller than $M_{h,k,m;\mu}$ at $x_0\in I$, then
\Eq{2+}{
  \frac{\partial_1^2\D_{f,g}(x_0,x_0)}{\partial_1\D_{f,g}(x_0,x_0)}
  \leq \frac{\partial_1^2\D_{h,k}(x_0,x_0)}{\partial_1\D_{h,k}(x_0,x_0)}.
}
On the other hand, if the functions
\Eq{*}{
  t\mapsto \partial_i m(\pmb{x}_0,t)-\langle \partial_im,1\rangle_\mu(\pmb{x}_0) \qquad(i\in\{1,\dots,d-1\})
}
are $\mu$-linearly independent and \eq{2+} holds with strict inequality, then $M_{f,g,m;\mu}$ is locally
smaller than $M_{h,k,m;\mu}$ at $x_0\in I$. }

\begin{proof}Assume that $M_{f,g,m;\mu}$ is locally smaller than $M_{h,k,m;\mu}$ at $x_0\in I$.
Then, by \cor{1}, the $(d-1)\times(d-1)$-matrix whose $(i,j)$th entry is given by
\Eq{8}{
  &\big(\langle \partial_im,\partial_j m\rangle_\mu
    -\langle\partial_im,1\rangle_\mu\langle\partial_jm,1\rangle_\mu\big)(\pmb{x}_0)\cdot
    \bigg(\frac{\partial_1^2\D_{h,k}(x_0,x_0)}{\partial_1\D_{h,k}(x_0,x_0)}
       - \frac{\partial_1^2\D_{f,g}(x_0,x_0)}{\partial_1\D_{f,g}(x_0,x_0)}\bigg)
}
for $i,j\in\{1,\dots,d-1\}$ is positive semidefinite at $x_0\in I$. This implies that all the diagonal
elements of this matrix are nonnegative, i.e., for all $i\in\{1,\dots,d-1\}$,
\Eq{8i}{
  &\big(\langle \partial_im,\partial_i m\rangle_\mu
    -\langle\partial_im,1\rangle_\mu\langle\partial_im,1\rangle_\mu\big)(\pmb{x}_0)\cdot
    \bigg(\frac{\partial_1^2\D_{h,k}(x_0,x_0)}{\partial_1\D_{h,k}(x_0,x_0)}
       - \frac{\partial_1^2\D_{f,g}(x_0,x_0)}{\partial_1\D_{f,g}(x_0,x_0)}\bigg)\geq0.
}
If, for some $i\in\{1,\dots,d-1\}$, the map $t\mapsto\partial_im(\pmb{x_0},t)$ is not $\mu$-almost
everywhere constant, then
\Eq{*}{
\mu(\{t\mid\partial_im(\pmb{x}_0,t)\neq\langle\partial_im,1\rangle_\mu(\pmb{x}_0)\})>0,
}
whence
\Eq{ii}{
\big(\langle \partial_im,\partial_i m\rangle_\mu
    -\langle\partial_im,1\rangle_\mu\langle\partial_im,1\rangle_\mu\big)(\pmb{x}_0)
=\int_T\big(\partial_im(\pmb{x_0},t)-\langle\partial_im,1\rangle_\mu(\pmb{x}_0)\big)^2\d\mu(t)>0.
}
Inequality \eq{ii}, combined with \eq{8i}, implies that
\Eq{*}{
\frac{\partial_1^2\D_{h,k}(x_0,x_0)}{\partial_1\D_{h,k}(x_0,x_0)}
       - \frac{\partial_1^2\D_{f,g}(x_0,x_0)}{\partial_1\D_{f,g}(x_0,x_0)}\geq0,
}
i.e, the inequality \eq{2+} holds.

Now assume that the functions
\Eq{fn}{
  t\mapsto \partial_i m(\pmb{x}_0,t)-\langle \partial_im,1\rangle_\mu(\pmb{x}_0) \qquad(i\in\{1,\dots,d-1\})
}
are $\mu$-linearly independent and \eq{2+} holds with strict inequality. It is clear that the $(d-1)\times(d-1)$-matrix
whose $(i,j)$th entry is given by
\Eq{Gr}{
\big(\langle \partial_im,\partial_j m\rangle_\mu
    -&\langle\partial_im,1\rangle_\mu\langle\partial_jm,1\rangle_\mu\big)(\pmb{x}_0)\\
    &=\int_T\big(\partial_im(\pmb{x_0},t)-\langle\partial_im,1\rangle_\mu(\pmb{x}_0)\big)
    \big(\partial_jm(\pmb{x_0},t)-\langle\partial_jm,1\rangle_\mu(\pmb{x}_0)\big)\d\mu(t)
}
for $i,j\in\{1,\dots,d-1\}$ is a so-called Gram matrix which is always positive semmidefinite (see
\cite{Bha07}). Since the functions \eq{fn} are $\mu$-linearly independent it follows that the Gram matrix
with entries given by \eq{Gr} is positive definite. This result, combined with the strict
inequality \eq{2+}, implies that the $(d-1)\times(d-1)$-matrix whose $(i,j)$th entry is given by \eq{8}
is positive definite at $x_0\in I$. Hence, by \cor{1}, $M_{f,g,m;\mu}$ is locally smaller than
$M_{h,k,m;\mu}$ at $x_0\in I$.
\end{proof}

Now we formulate a particular case concerning generalized Gini means when the partial derivatives can be
calculated more explicitly. Indeed, if, for given $p,q\in\R$, the functions $f$ and $g$ are given by 
equations \eq{pq}, then
\Eq{De}{
   \D_{f,g}(x,y)=\Delta_{p,q}(x,y):=y^{p+q}\delta_{p,q}\Big(\frac{x}{y}\Big)
      \qquad(x,y\in\R_+),
}
where
\Eq{de}{
   \delta_{p,q}(t)
     :=\begin{cases}
        \dfrac{t^p-t^q}{p-q} &\mbox{if } p\neq q\\[3mm]
        t^p\ln t             &\mbox{if } p=q
      \end{cases}
      \qquad(t\in\R_+).
}
Then, one can easily get that $\delta'_{p,q}(1)=1$ and $\delta''_{p,q}(1)=p+q-1$, whence 
\Eq{*}{
  \partial_1\Delta_{p,q}(x,x)=x^{p+q-1} \qquad\mbox{and}\qquad 
  \partial_1^2\Delta_{p,q}(x,x)=(p+q-1)x^{p+q-2}.
}
Therefore,
\Eq{DD}{
   \frac{\partial_1^2\D_{f,g}(x,x)}{\partial_1\D_{f,g}(x,x)}
   =\frac{\partial_1^2\Delta_{p,q}(x,x)}{\partial_1\Delta_{p,q}(x,x)}
    =\frac{(p+q-1)x^{p+q-2}}{x^{p+q-1}}=(p+q-1)\frac{1}{x}
    \qquad(x\in\R_+).
}

\Cor{G1}{Assume that $I\subseteq \R_+$. Let $p,q,r,s\in\R$, let $m\in\C_2(I^d\times T)$ 
and $n\in\C_2(I^d\times S)$ be measurable families of means, and let $\mu$ and $\nu$ be probability measures 
on the measurable spaces $(T,\A)$ and $(S,\B)$, respectively. Suppose that $G_{p,q,m;\mu}$ is locally smaller 
than $G_{r,s,n;\nu}$ at $x_0\in I$. Then, there exists a neighborhood $U\subseteq I$ of $x_0$ such that for 
$x\in U$, for all $i\in\{1,\dots,d\}$, \eq{6} holds. In addition, the
$(d-1)\times(d-1)$-matrix whose $(i,j)$th entry is given by
\Eq{G7}{
  &\big(\langle \partial_in,\partial_j n\rangle_\nu
    -\langle\partial_in,1\rangle_\nu\langle\partial_jn,1\rangle_\nu\big)(\pmb{x}_0)
    (p+q-1)\frac{1}{x_0}
    +\langle \partial_i\partial_j n,1\rangle_\nu(\pmb{x}_0)\\
    &-\big(\langle \partial_im,\partial_j m\rangle_\mu
    -\langle\partial_im,1\rangle_\mu\langle\partial_jm,1\rangle_\mu\big)(\pmb{x}_0)
    (r+s-1)\frac{1}{x_0}
    -\langle \partial_i\partial_j m,1\rangle_\mu(\pmb{x}_0)
}
for $i,j\in\{1,\dots,d-1\}$ is positive semidefinite. \\\indent On the other hand, if \eq{6} holds for all 
$i\in\{1,\dots,d\}$ and for all $x$ in a neighborhood of $x_0$ and the $(d-1)\times (d-1)$-matrix whose 
$(i,j)$th entry is given by \eq{G7} is positive definite, then $G_{p,q,m;\mu}$ is locally smaller than 
$G_{r,s,n;\nu}$ at $x_0\in I$. }

\begin{proof}
The proof is a direct consequence of \cor{1} and formula \eq{DD}.
\end{proof}

\Cor{G2}{Assume that $I\subseteq \R_+$. Let $p,q,r,s\in\R$, let $m\in\C_2(I^d\times T)$ be a measurable
family of means, and let $\mu$ be a probability measure on the measurable space $(T,\A)$. Let $x_0\in I$ and
assume that there exists $i\in\{1,\dots,d-1\}$ such that, the map $t\mapsto\partial_im(\pmb{x_0},t)$ is not
$\mu$-almost everywhere constant on $T$. If $G_{p,q,m;\mu}$ is locally smaller than $G_{r,s,m;\mu}$ at
$x_0\in I$, then
\Eq{G2+}{
  p+q\leq r+s.
}
On the other hand, if $x_0\in I$, the functions
\Eq{*}{
  t\mapsto \partial_i m(\pmb{x}_0,t)-\langle \partial_im,1\rangle_\mu(\pmb{x}_0) \qquad(i\in\{1,\dots,d-1\})
}
are $\mu$-linearly independent and \eq{G2+} holds with strict inequality, then $G_{p,q,m;\mu}$ is locally
smaller than $G_{r,s,m;\mu}$ at $x_0\in I$. }

\begin{proof}
Applying \cor{2} and using formula \eq{DD}, the result follows immediately.
\end{proof}

\section{Necessary and sufficient conditions for global comparison of means}

In the rest of the paper, we consider the case when $\mu=\nu$ and $m=n$. In what follows, we give a condition 
containing two independent variables for \eq{2} which does not involve the measure $\mu$ and assumes 
first-order continuous differentiability of the Chebyshev system. In the special setting when $T=[0,1]$, 
$d=2$, $m$ is given by $m(\pmb{x},t):=tx_1+(1-t)x_2$, the following theorem simplifies to the result of 
\cite[Theorem 6]{LosPal08}. 

\Thm{2}{Let $(f,g),(h,k)\in\C_1(I)$ be Chebyshev systems, let $T$ be a compact and connected topological 
space and let $m:I^d\times T\to\R$ be a continuous family of $d$-variable means. Define the set $U_m$ by
\Eq{Um}{
  U_m:=\{(u,v)\mid \exists\, \pmb{x}\in I^d : u,v\in[\underline{m}(\pmb{x}),\overline{m}(\pmb{x})]\}
  =\bigcup_{\pmb{x}\in I^d}[\underline{m}(\pmb{x}),\overline{m}(\pmb{x})]^2.
}
The following three assertions are equivalent:
\begin{enumerate}[(i)]
\item for all Borel probability measures $\mu$ on $T$,
\Eq{13}{
  M_{f,g,m;\mu}(\pmb{x})\le M_{h,k,m;\mu}(\pmb{x})\qquad(\pmb{x}\in I^d);
}
\item there exists a nullsequence $(\gamma_j)$ of positive numbers in $[0,1]$ such that, for all $t_0,t\in T$
and for all $j\in\N$,
\Eq{13a}{
  M_{f,g,m;(1-\gamma_j)\delta_{t_0}+\gamma_j\delta_t}(\pmb{x})
   \le M_{h,k,m;(1-\gamma_j)\delta_{t_0}+\gamma_j\delta_t}(\pmb{x})\qquad(\pmb{x}\in I^d);
}
\item for all $(u,v)\in U_m$,
\Eq{12}{
 \frac{\D_{f,g}(u,v)}{\partial_1\D_{f,g}(v,v)}
  \le \frac{\D_{h,k}(u,v)}{\partial_1\D_{h,k}(v,v)}.
}
\end{enumerate}}

\begin{proof}
The implication (i)$\Longrightarrow$(ii) is obvious.

To prove (ii)$\Longrightarrow$(iii), let $(u,v)\in U_m$. Then there exists $\pmb{x}\in I^d$
such that $u,v\in[\underline{m}(\pmb{x}),\overline{m}(\pmb{x})]$. Due to the compactness and connectedness of
$T$, we have that \eq{mm} holds.
Therefore, there exits $t_0,t\in T$ such that $u=m(\pmb{x},t)$ and $v=m(\pmb{x},t_0)$. Applying \lem{2} twice
with the measure sequence  $\mu_j:=(1-\gamma_j)\delta_{t_0}+\gamma_j\delta_t$ and using inequality \eq{13}, we get
\Eq{*}{
 \frac{\D_{f,g}(u,v)}{\partial_1\D_{f,g}(v,v)}
   &=\frac{\D_{f,g}(m(\pmb{x},t),m(\pmb{x},t_0))}{\partial_1\D_{f,g}(m(\pmb{x},t_0),m(\pmb{x},t_0))}
   =\lim_{j\to\infty}
    \frac{1}{\gamma_j}\left[M_{f,g,m;\mu_j}(\pmb{x})-m(\pmb{x},t_0)\right]\\
   &\le\lim_{j\to\infty}
    \frac{1}{\gamma_j}\left[M_{h,k,m;\mu_j}(\pmb{x})-m(\pmb{x},t_0)\right]
   =\frac{\D_{h,k}(m(\pmb{x},t),m(\pmb{x},t_0))}{\partial_1\D_{h,k}(m(\pmb{x},t_0),m(\pmb{x},t_0))}
   =\frac{\D_{h,k}(u,v)}{\partial_1\D_{h,k}(v,v)},
}
which proves \eq{12}.

For the proof of (iii)$\Longrightarrow$(i), let $\pmb{x}\in I^d$ be arbitrarily fixed. In view of the
inclusion $M_{h,k,m;\mu}(\pmb{x})\in[\underline{m}(\pmb{x}),\overline{m}(\pmb{x})]$ and the equality \eq{mm},
there exits $t_0\in T$ such that
\Eq{*}{
  m(\pmb{x},t_0):=M_{h,k,m;\mu}(\pmb{x}).
}
Taking any $t\in T$ and applying inequality \eq{12} for $u:=m(\pmb{x},t)$ and $v:=m(\pmb{x},t_0)$, we get
\Eq{*}{
 \frac{\D_{f,g}(m(\pmb{x},t),m(\pmb{x},t_0))}{\partial_1\D_{f,g}(m(\pmb{x},t_0),m(\pmb{x},t_0))}
  \le \frac{\D_{h,k}(m(\pmb{x},t),m(\pmb{x},t_0))}{\partial_1\D_{h,k}(m(\pmb{x},t_0),m(\pmb{x},t_0))}.
}
Integrating this inequality with respect to the variable $t\in T$, we get
\Eq{15}{
  \dfrac{\int_T\D_{f,g}\big(m(\pmb{x},t),M_{h,k,m;\mu}(\pmb{x})\big)\d\mu(t)}
    {\partial_1\D_{f,g}\big(M_{h,k,m;\mu}(\pmb{x}),M_{h,k,m;\mu}(\pmb{x})\big)}
  \le\dfrac{\int_T\D_{h,k}\big(m(\pmb{x},t),M_{h,k,m;\mu}(\pmb{x})\big)\d\mu(t)}
    {\partial_1\D_{h,k}\big(M_{h,k,m;\mu}(\pmb{x}),M_{h,k,m;\mu}(\pmb{x})\big)}.
}
By the definition of the value $M_{h,k,m;\mu}(\pmb{x})$, the numerator of the right hand side of this
inequality is zero, whence we obtain
\Eq{*}{
\dfrac{\int_T\D_{f,g}\big(m(\pmb{x},t),M_{h,k,m;\mu}(\pmb{x})\big)\d\mu(t)}
    {\partial_1\D_{f,g}\big(M_{h,k,m;\mu}(\pmb{x}),M_{h,k,m;\mu}(\pmb{x})\big)}\leq0.
}
If $\partial_1\D_{f,g}(x,x)<0$ for all $x\in I$, then
\Eq{last}{
\int_T\D_{f,g}\big(m(\pmb{x},t),M_{h,k,m;\mu}(\pmb{x})\big)\d\mu(t)\geq0
}
and also $(f,g)$ is a positive Chebyshev system, hence, by \lem{0}, the above inequality implies that
\eq{13} holds. In the other possible case, i.e., when $\partial_1\D_{f,g}(x,x)>0$ for all $x\in I$, then
inequality \eq{last} is reversed but $(f,g)$ is a negative Chebyshev system, thus by \lem{0} again, inequality
\eq{13} follows.
\end{proof}

Having a look at the proof, one can see that the compactness and connectedness of $T$ was only used to prove 
implication (ii)$\Longrightarrow$(iii).

\Cor{pq1}{Assume that $I\subseteq \R_+$ and $p,q,r,s\in\R$. Let $T$ be a compact and connected topological 
space and let $m:I^d\times T\to\R$ be a continuous family of means. Define the constant $m^*\in[1,+\infty]$ by
\Eq{*}{
  m^*:=\sup_{\pmb{x}\in I^d}\,\frac{\overline{m}(\pmb{x})}{\underline{m}(\pmb{x})}.
}
The following three assertions are equivalent:
\begin{enumerate}[(i)]
\item for all Borel probability measures $\mu$ on $T$,
\Eq{G13}{
  G_{p,q,m;\mu}(\pmb{x})\le G_{r,s,m;\mu}(\pmb{x})\qquad(\pmb{x}\in I^d);
}
\item there exists a null sequence $(\gamma_j)$ of positive numbers in $[0,1]$ such that, for all $t_0,t\in T$
and for all $j\in\N$,
\Eq{G13a}{
  G_{p,q,m;(1-\gamma_j)\delta_{t_0}+\gamma_j\delta_t}(\pmb{x})
   \le G_{r,s,m;(1-\gamma_j)\delta_{t_0}+\gamma_j\delta_t}(\pmb{x})\qquad(\pmb{x}\in I^d);
}
\item
\Eq{14}{
  \delta_{p,q}(t)\leq\delta_{r,s}(t) \qquad\big(t\in\,\big](m^*)^{-1},m^*\big[\big);
}
\item In the case $m^*=+\infty$,
\Eq{G12}{
  \min(p,q)\leq\min(r,s)\qquad\mbox{and}\qquad \max(p,q)\leq\max(r,s),
}
while in the case $m^*<+\infty$,
\Eq{G15}{
  \delta_{p,q}\big((m^*)^{-1}\big)\leq\delta_{r,s}\big((m^*)^{-1}\big),\qquad
  \delta_{p,q}\big(m^*\big)\leq\delta_{r,s}\big(m^*\big),
  \qquad\mbox{and}\qquad p+q\leq r+s.
}
\end{enumerate}}

\begin{proof}
Applying \thm{2} and using notations introduced in \eq{De} and \eq{de} imply that conditions \eq{G13} and 
\eq{G13a} are equivalent to the inequality
\Eq{*}{
   \frac{\Delta_{p,q}(v,u)}{\partial_1\Delta_{p,q}(u,u)}
   \leq \frac{\Delta_{r,s}(v,u)}{\partial_1\Delta_{r,s}(u,u)}
   \qquad((u,v)\in U_m),
}
where the set $U_m$ is defined in \eq{Um}. This inequality can be rewritten as
\Eq{de+}{
  v\delta_{p,q}\Big(\frac{u}{v}\Big)\leq v\delta_{r,s}\Big(\frac{u}{v}\Big) \qquad((u,v)\in U_m).
}
Observe that
\Eq{mm+}{
  \big](m^*)^{-1},m^*\big[\,\subseteq \Big\{\frac{u}{v}\colon (u,v)\in U_m\Big\}
  \subseteq\,\big[(m^*)^{-1},m^*\big].
}
Indeed, if $t\in \big](m^*)^{-1},m^*\big[$ and $t\geq1$, then $t<m^*$, hence there exits $\pmb{x}\in I^d$ 
such that $t<\frac{\overline{m}(\pmb{x})}{\underline{m}(\pmb{x})}$. Then, with $v=\underline{m}(\pmb{x})$, 
$u=t\underline{m}(\pmb{x})$, we have that $t=\frac{u}{v}$ and 
$u,v\in[\underline{m}(\pmb{x}),\overline{m}(\pmb{x})]$. Therefore $t$ is of the form $\frac{u}{v}$ for some 
$(u,v)\in U_m$. A similar argument yields for $t\leq1$ a similar representation. This proves the first 
inclusion.  

For the second inclusion, observe that if $(u,v)\in U_m$, then, for some $\pmb{x}\in I^d$, we have
$u,v\in[\underline{m}(\pmb{x}),\overline{m}(\pmb{x})]$. Hence
\Eq{*}{
  (m^*)^{-1}\leq\frac{\underline{m}(\pmb{x})}{\overline{m}(\pmb{x})}
  \leq\frac{u}{v}\leq \frac{\overline{m}(\pmb{x})}{\underline{m}(\pmb{x})}\leq m^*
}
Therefore, in view of the inclusions in \eq{mm+}, inequality \eq{de} is equivalent to condition \eq{14}. 

To show the equivalence of condition (iv) to the previous ones, we have to distinguish two cases.
If $m^*=+\infty$, then $(m^*)^{-1}=0$, therefore (iii) can be rewritten as
\Eq{*}{
  \delta_{p,q}(t)\leq\delta_{r,s}(t) \qquad\big(t\in\,]0,\infty[\big).
}
This inequality is known to be equivalent (cf.\ \cite{DarLos70}) to the comparison inequality
\Eq{*}{
  G_{p,q}(\pmb{x})\leq G_{r,s}(\pmb{x}) \qquad(d\in\N,\,\pmb{x}\in\big]0,\infty[^d),
}
of Gini means (with arbitrary many variables over the interval $]0,\infty[$). In view of the result \cite[Satz 
5]{DarLos70}, the above inequality is characterized by the condition \eq{G12}. Therefore (iii) is equivalent 
to (iv) in this case.

Now consider the case $m^*<+\infty$. Then the inequality in (iii) is equivalent to the comparison inequality
\Eq{*}{
  G_{p,q}(\pmb{x})\leq G_{r,s}(\pmb{x}) \qquad(d\in\N,\,\pmb{x}\in\big]1,m^*[^d),
}
of Gini means (with arbitrary many variables over the interval $]1,m^*[$). Using the results of the papers 
\cite[Theorem 7]{Los77} or \cite{Pal89c}, it follows that the above inequality is characterized by \eq{G15}, 
which implies that (iii) is equivalent to (iv) also in this case.

This completes the proof of the corollary.
\end{proof}

\section{Necessary and sufficient conditions for the local and global comparison of generalized 
quasi-arithmetic means}

In the next result we offer 6 equivalent conditions for the comparison of $d$-variable generalized 
quasi-arithmetic means. The interesting feature of this result is the equivalence of the global and local 
comparability.

\Thm{Q}{Let $f,h:I\to\R$ be twice continuously differentiable functions with non-vanishing first derivatives, 
and let $m\in\C_2(I^d\times T)$ be a measurable family of $d$-variable means. Let $\mu_0$ be a probability 
measure such that, for all $x_0\in I$, there exists $i\in\{1,\dots,d-1\}$ such that 
$t\mapsto\partial_im(\pmb{x}_0,t)$ is  not $\mu_0$-almost everywhere constant on $T$. 
The following assertions are equivalent:
\begin{enumerate}[(i)]
\item for all Borel probability measures $\mu$ on $T$,
\Eq{Q1}{
  M_{f,1,m;\mu}(\pmb{x})\le M_{h,1,m;\mu}(\pmb{x})\qquad(\pmb{x}\in I^d);
}
\item 
\Eq{*}{
  M_{f,1,m;\mu_0}(\pmb{x})
   \le M_{h,1,m;\mu_0}(\pmb{x})\qquad(\pmb{x}\in I^d);
}
\item for all $x_0\in I$, there exists a neighborhood $U\subseteq I$ of $x_0$ such that
\Eq{*}{
  M_{f,1,m;\mu_0}(\pmb{x})
   \le M_{h,1,m;\mu_0}(\pmb{x})\qquad(\pmb{x}\in U^d);
}
\item for all $x\in I$,
\Eq{Q4}{
  \frac{f''(x)}{f'(x)}\leq \frac{h''(x)}{h'(x)};
}
\item the function $h\circ f^{-1}$ is convex (concave) on $f(I)$ provided that $f$ is increasing (decreasing);
\item for all $(u,v)\in I^2$,
\Eq{Q5}{
  \frac{f(u)-f(v)}{f'(v)} \leq \frac{h(u)-h(v)}{h'(v)}.
}
\end{enumerate}}

\begin{proof} The implications (i)$\Longrightarrow$(ii) and (ii)$\Longrightarrow$(iii) are trivial.

To prove (iii)$\Longrightarrow$(iv), will apply \cor{2}. Let $x_0\in I$ be arbitrary. Then (iii) asserts that 
$M_{f,1,m;\mu_0}$ is locally smaller than $M_{h,1,m;\mu_0}$ at $x_0$ and we also have an index 
$i\in\{1,\dots,d-1\}$ such that $t\mapsto\partial_im(\pmb{x}_0,t)$ is  not $\mu_0$-almost everywhere constant 
on $T$. Therefore, by \cor{2}, inequality \eq{2+} follows with the functions $g:=k:=1$. It is immediate to 
see that \eq{2+} implies \eq{Q4} for $x=x_0$.

Now assume (iv) and that $f$ is increasing (the nondecreasing case is analogous). Then $g:=h\circ f^{-1}$ is 
twice differentiable on $f(I)$. By \eq{Q4}, the ratio $\frac{h'}{f'}$ is a nondecreasing function. Therefore, 
$g'=\frac{h'\circ f^{-1}}{f'\circ f^{-1}}$ is also nondecreasing, which proves that $g''\geq0$. Hence $g$ 
must be convex on $f(I)$, i.e., (v) holds.

If (v) is valid and $f$ is increasing, then the convexity of $g:=h\circ f^{-1}$ implies that
\Eq{*}{
   g(y)+g'(y)(x-y)\leq g(x) 
}
for all $x,y\in f(I)$. With the substitution $x=f(u)$, $y=f(v)$, where $u,v\in I$, the above inequality 
reduces to \eq{Q5}, proving (vi).

Finally, assume that (vi) holds. Observe that then \eq{12} is valid for all $(u,v)\in U_m$ with the 
functions $g:=k:=1$. Thus, the condition (iii) of \thm{2} is satisfied, whence it follows that the mean
$M_{f,1,m;\mu_0}$ is (globally) smaller than $M_{h,1,m;\mu_0}$ on $I^d$, i.e., \eq{Q1} holds.
\end{proof}

As an immediate consequence, we obtain the characterization of the comparison among generalized H\"older means.

\Cor{Q+}{Let $I\subseteq\R_+$, $p,q\in\R$, and let $m\in\C_2(I^d\times T)$ be a measurable family of 
$d$-variable means. Let $\mu_0$ be a probability measure such that, for all $x_0\in I$, there exists 
$i\in\{1,\dots,d-1\}$ such that $t\mapsto\partial_im(\pmb{x}_0,t)$ is  not $\mu_0$-almost everywhere constant 
on $T$. The following assertions are equivalent:
\begin{enumerate}[(i)]
\item for all Borel probability measures $\mu$ on $T$,
\Eq{*}{
  G_{p,0,m;\mu}(\pmb{x})\le G_{q,0,m;\mu}(\pmb{x})\qquad(\pmb{x}\in I^d);
}
\item 
\Eq{*}{
  G_{p,0,m;\mu_0}(\pmb{x})
   \le G_{q,0,m;\mu_0}(\pmb{x})\qquad(\pmb{x}\in I^d);
}
\item for all $x_0\in I$, there exists a neighborhood $U\subseteq I$ of $x_0$ such that
\Eq{*}{
  G_{p,0,m;\mu_0}(\pmb{x})
   \le G_{q,0,m;\mu_0}(\pmb{x})\qquad(\pmb{x}\in U^d);
}
\item $p\leq q$.
\end{enumerate}}

\begin{proof}
By taking $f(x):=x^p$ if $p\neq0$ and $f(x):=\log(x)$ if $p=0$ and $h(x):=x^q$ if $q\neq0$ and 
$h(x):=\log(x)$ if $q=0$ and applying \thm{Q} the result follows immediately because conditions (i), (ii) and 
(iii) are equivalent to the same conditions of \thm{Q}, and $p\leq q$ is equivalent to condition (iv) of 
\thm{Q}.
\end{proof}

\def\MR{}

\begin{thebibliography}{10}

\bibitem{Baj58}
M.~Bajraktarević.
\newblock {Sur une équation fonctionnelle aux valeurs moyennes}.
\newblock {\em Glasnik Mat.-Fiz. Astronom. Društvo Mat. Fiz. Hrvatske Ser.
  II}, 13:243–248, 1958.

\bibitem{Baj69}
M.~Bajraktarević.
\newblock {Über die {V}ergleichbarkeit der mit {G}ewichtsfunktionen gebildeten
  {M}ittelwerte}.
\newblock {\em Studia Sci. Math. Hungar.}, 4:3–8, 1969.

\bibitem{Ber98}
L.~R. Berrone.
\newblock {The mean value theorem: functional equations and {L}agrangian
  means}.
\newblock {\em Epsilon}, 14(1(40)):131–151, 1998.

\bibitem{BerMor98}
L.~R. Berrone and J.~Moro.
\newblock {Lagrangian means}.
\newblock {\em Aequationes Math.}, 55(3):217–226, 1998.

\bibitem{BesPal03}
M.~Bessenyei and Zs. Páles.
\newblock {Hadamard-type inequalities for generalized convex functions}.
\newblock {\em Math. Inequal. Appl.}, 6(3):379–392, 2003.

\bibitem{Bha07}
R.~Bhatia.
\newblock {\em Positive definite matrices}.
\newblock Princeton Series in Applied Mathematics. Princeton University Press,
  Princeton, NJ, 2007.

\bibitem{CziPal04}
P.~Czinder and Zs. Páles.
\newblock {An extension of the {H}ermite-{H}adamard inequality and an
  application for {G}ini and {S}tolarsky means}.
\newblock {\em J. Inequal. Pure Appl. Math.}, 5(2):Article 42, pp. 8
  (electronic), 2004.

\bibitem{CziPal05}
P.~Czinder and Zs. Páles.
\newblock {Local monotonicity properties of two-variable {G}ini means and the
  comparison theorem revisited}.
\newblock {\em J. Math. Anal. Appl.}, 301(2):427–438, 2005.

\bibitem{CziPal06}
P.~Czinder and Zs. Páles.
\newblock {Some comparison inequalities for {G}ini and {S}tolarsky means}.
\newblock {\em Math. Inequal. Appl.}, 9(4):607–616, 2006.

\bibitem{DarLos70}
Z.~Daróczy and L.~Losonczi.
\newblock {Über den {V}ergleich von {M}ittelwerten}.
\newblock {\em Publ. Math. Debrecen}, 17:289–297 (1971), 1970.

\bibitem{Gin38}
C.~Gini.
\newblock {{D}i una formula compressiva delle medie}.
\newblock {\em Metron}, 13:3–22, 1938.

\bibitem{HarLitPol34}
G.~H. Hardy, J.~E. Littlewood, and G.~Pólya.
\newblock {\em {Inequalities}}.
\newblock Cambridge University Press, Cambridge, 1934.
\newblock (first edition), 1952 (second edition).

\bibitem{LeaSho78}
E.~Leach and M.~Sholander.
\newblock {Extended mean values}.
\newblock {\em Amer. Math. Monthly}, 85(2):84–90, 1978.

\bibitem{LeaSho83}
E.~Leach and M.~Sholander.
\newblock {Extended mean values {I}{I}}.
\newblock {\em J. Math. Anal. Appl.}, 92:207–223, 1983.

\bibitem{Los77}
L.~Losonczi.
\newblock {Inequalities for integral mean values}.
\newblock {\em J. Math. Anal. Appl.}, 61(3):586–606, 1977.

\bibitem{Los99}
L.~Losonczi.
\newblock {Equality of two variable weighted means: reduction to differential
  equations}.
\newblock {\em Aequationes Math.}, 58(3):223–241, 1999.

\bibitem{Los00a}
L.~Losonczi.
\newblock {Equality of {C}auchy mean values}.
\newblock {\em Publ. Math. Debrecen}, 57:217–230, 2000.

\bibitem{Los02d}
L.~Losonczi.
\newblock {Comparison and subhomogeneity of integral means}.
\newblock {\em Math. Inequal. Appl.}, 5(4):609–618, 2002.

\bibitem{Los02c}
L.~Losonczi.
\newblock {On the comparison of {C}auchy mean values}.
\newblock {\em J. Inequal. Appl.}, 7(1):11–24, 2002.

\bibitem{LosPal08}
L.~Losonczi and Zs. Páles.
\newblock {Comparison of means generated by two functions and a measure}.
\newblock {\em J. Math. Anal. Appl.}, 345(1):135–146, 2008.

\bibitem{NeuPal03}
E.~Neuman and Zs. Páles.
\newblock {On comparison of {S}tolarsky and {G}ini means}.
\newblock {\em J. Math. Anal. Appl.}, 278(2):274–284, 2003.

\bibitem{Pal88b}
Zs. Páles.
\newblock {Inequalities for differences of powers}.
\newblock {\em J. Math. Anal. Appl.}, 131(1):271–281, 1988.

\bibitem{Pal88c}
Zs. Páles.
\newblock {Inequalities for sums of powers}.
\newblock {\em J. Math. Anal. Appl.}, 131(1):265–270, 1988.

\bibitem{Pal89c}
Zs. Páles.
\newblock {On comparison of homogeneous means}.
\newblock {\em Ann. Univ. Sci. Budapest. Eötvös Sect. Math.}, 32:261–266
  (1990), 1989.

\bibitem{Pal92a}
Zs. Páles.
\newblock {Comparison of two variable homogeneous means}.
\newblock In W.~Walter, editor, {\em {General Inequalities, 6 (Oberwolfach,
  1990)}}, {International Series of Numerical Mathematics}, page 59–70.
  Birkhäuser, Basel, 1992.

\bibitem{Sto75}
K.~B. Stolarsky.
\newblock {{G}eneralizations of the logarithmic mean}.
\newblock {\em Math. Mag.}, 48:87–92, 1975.

\end{thebibliography}

\end{document}